\documentclass[review]{siamltex}
\usepackage{latexsym, graphicx, epsfig, amsmath, amsfonts,amssymb,subfigure}
\usepackage{multirow}
\usepackage{color}
 \usepackage[percent]{overpic}

\usepackage{algorithm}
\usepackage{algorithmicx}
\usepackage{algpseudocode}
\floatname{algorithm}{Algorithm}
\usepackage{threeparttable}
\usepackage{booktabs}

\def\mb{\mathbf}
\def\R{\mathbb{R}}

\title{An adaptive multi-fidelity PC-based ensemble Kalman inversion  for inverse problems}%

\author{Liang Yan\thanks{Department of Mathematics,  Southeast University, Nanjing, 210096, China (yanliang@seu.edu.cn). This author's work was supported by NSF of China (No.11771081) and Qing Lan project of Jiangsu Province. }
          \and Tao Zhou\thanks{LSEC, Institute of Computational Mathematics and Scientific/Engineering Computing, Academy of Mathematics and Systems Science, Chinese Academy of Sciences, Beijing 100190, China (tzhou@lsec.cc.ac.cn). This author's work is partially supported by the NSF of China (under grant numbers 11822111, 11688101, 91630203, 11571351, and 11731006), the science challenge project (No. TZ2018001), the national key basic research program (No. 2018YFB0704304), NCMIS, and the youth innovation promotion association (CAS) }}
\begin{document}
%\graphicspath{figure/}
\maketitle

\begin{abstract}
The ensemble Kalman inversion (EKI), as a derivative-free methodology, has been widely used  in the parameter estimation of inverse problems. Unfortunately, its cost may become moderately large for systems described by high dimensional nonlinear PDEs, as EKI requires a relatively large ensemble size to guarantee its performance. In this paper, we propose an adaptive multi-fidelity polynomial chaos (PC)  based EKI technique to address  this challenge. Our new strategy  combines a large number of low-order PC surrogate model evaluations and a small number of high-fidelity forward model evaluations, yielding a multi-fidelity approach. Especially, we present a new approach that adaptively constructs and refines a multi-fidelity PC surrogate  during the EKI simulation.  Since the forward model evaluations are only required for updating the low-order multi-fidelity PC model, whose number can be much smaller than the total ensemble size of the classic EKI, the entire computational costs are thus significantly reduced. The new algorithm was tested through the two-dimensional time fractional  inverse diffusion problems and demonstrated great effectiveness in comparison with PC based EKI and classic EKI.
\end{abstract}


\pagestyle{myheadings}
\thispagestyle{plain}
\markboth{LIANG YAN AND TAO ZHOU}
{ AMPC-EKI FOR INVERSE PROBLES}

\section{Introduction}

The estimation of model parameters from a set of observations is a key requirement in science and engineering. In practical applications, the observations are always indirect, noisy, and limited in number. Quantifying the resulting uncertainty in parameters is then an indispensable part of the inference process. While parameter estimation problems can be solved using a variety of approaches \cite{Tarantola2005inverse}, the Bayesian approach \cite{Kaipio+Somersalo2005,Stuart2010} is particularly attractive, as it provides a systematic framework for quantifying parameter uncertainty. Moreover, the Bayesian approach can incorporate uncertainties in the model and the observations and leads to a complete characterization of the uncertainty in terms of the posterior distribution, i.e. the conditional distribution of the unknown parameters given the data.  Since the posterior is typically not of analytical form and cannot be easily interrogated, many numerical approaches such as Markov chain Monte Carlo (MCMC) methods have been developed.   However, for computationally intensive applications, the computation of the posterior is prohibitively expensive and is even  intractable.

The ensemble Kalman filter (EnKF) \cite{Evensen1994}, which is a Monte Carlo variant of the classical Kalman filter \cite{kalman1960new}, is a computationally alternative for parameter estimation of inverse problems. Since it only requires the evaluation of the forward operator but not its derivative, this approach has been successfully used in various areas of application, including oceanic \cite{Bertino+Evensen2003sequential} and geophysical \cite{Aanonsen2009ensemble}.  As a sequential data assimilation technique, EnKF needs to update model parameters and states simultaneously at each assimilation step, which makes its application inconvenient when the model involves multiple processes \cite{Emerick+Reynolds2013}. In this situation, computing a global update with all available data is preferred, and this leads to the scheme of ensemble smoother (ES)\cite{Evensen2009book,Van+Evensen1996}. It has been shown that ES can obtain comparable results as EnKF in some parameter estimation problems \cite{Skjervheim+Evensen2011}.  However, for strongly nonlinear problems, both EnKF and ES \cite{Emerick+Reynolds2012, Emerick+Reynolds2013,Gu+Oliver2007,Lorentzen+Naevdal2011} need certain iterations to achieve satisfactory data matches.  Variants ideas including ensemble randomized likelihood \cite{Chen+Oliver2012ensemble,Chen+Oliver2013levenberg} and multiple-data-assimilation ES \cite{Emerick+Reynolds2012, Emerick+Reynolds2013} have been investigated along this line.  The iterations of the smoother update turn out to partly resolve issues with nonlinearity and lead to better results than what is obtained by ES \cite{Evensen2018analysis}. The methodology is described in a basic form, applicable to a general  inverse problem in \cite{Iglesias+Law+Kody2013ensemble}. It is widely known the lack of stability of ensemble Kalman-based methods when the ensemble size is small with respect to the number of parameters or measurements. Therefore, many attentions have been recently given to the regularization of ensemble Kalman-based methods by means of localization and inflation \cite{Asch2016data,Law+Stuart2015data}.  It is also important to mention the work of \cite{Iglesias2016regularizing} that uses a regularizing iterative ensemble Kalman method to solve PDE-constrained parameter identification problems.

While both aforementioned ensemble Kalman-based inversion methods (EKIs) are derivate-free approaches where the ensemble updates are based on simple computations involving covariances and crosscovariances computed directly from the ensemble realizations. However,  as a Monte Carlo method, EKI requires a sufficiently large ensemble size to guarantee reliable estimations.  This is undesirable in practice as each ensemble realization requires a solution of the forward model and can be time-consuming to compute for large-scale complex systems.   One approach to alleviate the computational burden is to use surrogate models, which are constructed to approximate the outputs of the forward model at low computational cost, see  \cite{Asher2015review,Li+Xiu2009generalized,Li+Lin+Zhang2014adaptive,ju2018adaptive,Saad2009characterization} and the references therein.   One of the most popular approaches is to expand the high-fidelity model  in a generalized polynomial chaos (gPC) basis and approximate its coefficients via either a Galerkin approach or a collocation approach \cite{xiu2010book}. When the PC surrogate is obtained, we can generate a large ensemble of realizations without incurring a notable computational cost. Therefore,  the majority of the computational cost in PC-based EKIs  is spent on  building the PC representation of the high-fidelity model. Usually, this cost is much smaller than running the simulations for all realizations in EKIs if the dimensionality of the problem, i.e., the number of random variables, is relatively low. However, PC-based EKIs  may lose its advantage over EKIs for relatively high-dimensional problems because the number of PC basis functions grows very fast as the dimensionality increases.  In this case, we can only use a small order of PC basis to construct the surrogate otherwise PC expansions become expensive. This introduces a large model error unless the forward model is well represented by a low-order PC expansion.  On the other hand, constructing a sufficiently accurate PC surrogate model over the support of the prior distribution may not be possible in many practical problems, especially when the data contain information beyond what is assumed in the prior \cite{lu2015JCP,Yan+Zhang2017IP}.  For these problems, we need a carefully designed method to balance the accuracy and efficiency of the PC-based EKIs.

%This is our follow up work of our previous work \cite{Yan+Zhou2018adaptive}. 
In this paper we shall propose an adaptive  multi-fidelity PC  EKI algorithm  to study inverse problems for parameter estimation. The main focus of the proposed method is to address two sources that could potentially affect the accuracy and efficiency of the PC-based EKIs, i.e.,   the high dimensionality of the parameter space as well as the intrinsic properties of the forward model.  In our previous work \cite{Yan+Zhou2018adaptive}, we have designed an adaptive multi-fidelity PC MCMC algorithm that samples the posterior and selects the nodes for constructing the multi-fidelity PC simultaneously.  In this study, inspired by the recent progress of data assimilation and the multi-fidelity PC, we develop a new adaptive multi-fidelity PC-based EKI for parameter estimation of inverse problems. We construct a multi-fidelity PC surrogate by combining the low-order PC surrogate evaluations and the forward model evaluations. We also propose a strategy that refines the surrogate adaptively during the EKI simulation.   Specifically, in order to address ill-posedness due to small ensemble size, we consider the regularization iterative ensemble Kalman smoother \cite{Iglesias2016regularizing}.  We remark that although we focus our attention on the version of EKIs presented in \cite{Iglesias2016regularizing}, our methods can be easily applied to other  ensemble Kalman-based methods with simple modifications.

The structure of the paper is as follows. In the next section, we shall review the formulation of regularization iterative ensemble Kalman smoother and the PC-surrogate  approach to  EKI. In section 3, we shall propose an adaptive multi-fidelity PC approach to EKI. In section 4, we use  a two dimensional  time-fractional inverse diffusion problem to demonstrate the accuracy and efficiency of the proposed method. We finally give some concluding remarks in Section 5.

\section {Background and problem setup}\label{sec:setup}

In this section, we first give a brief overview of the  Bayesian inverse problems. Then we will introduce a regularization iterative ensemble Kalman smoother  and PC-based  EKI.

\subsection{Bayesian inverse problems}

The most standard approach to quantify the uncertainty in parameters is the Bayesian framework. The aim is to merge uncertainties, both in prior knowledge and observational data, with the mathematical model. The  prior belief about the parameter $\theta\in \R^d$ is encoded in the prior probability distribution $\pi(\theta)$. The data $y\in \R^m$ and the parameter $\theta$ are related via the forward model (also known as the  parameter-to-observable map) $f$ by
\begin{equation}\label{ipeq}
y = f(\theta)+\xi
\end{equation}
where $\xi \in \R^m$ is the measurement error.  We assume that the error $\xi$ is a Gaussian random vector with mean zero and covariance matrix $\Gamma \in \R^{m\times m}$, i.e., $\xi \sim N(0, \Gamma)$. The likelihood of the measure data $y$ given a particular instance of $\theta$ is denoted by $\pi(y|\theta)$.  In the Bayesian framework,  the  distribution of the $\theta$  conditioned on the data $y$, i.e., the posterior distribution $\pi(\theta|y)$ follows the Bayes' rule,
\begin{eqnarray}\label{ppdf}
\pi(\theta|y) \propto \pi(y|\theta) \pi(\theta)\propto \exp(-\Phi(\theta;y))\pi(\theta),
\end{eqnarray}
where the potential $\Phi(\theta;y)$ is defined by
\begin{equation}\label{poteq}
\Phi(\theta;y)=\|\Gamma^{-1/2}(y-f(\theta))\|.
\end{equation}

Since the forward model $f$ is always nonlinear, the expression of the potential yields a posterior distribution that cannot be written in a closed form.  Standard sampling methods, e.g. MCMC, have been extensively used to sample such unknown posterior distributions. Unfortunately, this approach often requires a large number of repeated evaluations of the forward model $f$, which can be very expensive.  The full characterization of the posterior by means of sampling is therefore impractical. In this work, we consider the application of iterative EKIs to approximate the Bayesian posterior. The connection between Bayesian inversion and the iterative  EKIs can be found in \cite{Evensen2018analysis,Iglesias2015iterative}.

\subsection{Regularization iterative ensemble Kalman smoother}
We follow closely  the framework in \cite{Iglesias2016regularizing}.  Assume we derive $N_e$ initial ensemble $\theta_0^{(j)}\, (j \in \{1, \cdots, N_e\})$ from the prior $\pi(\theta)$. When applied to the inverse problem (\ref{ipeq}), we use  the following iterative procedure, in which the subscript $n$ denotes the iteration step, and the superscript $(j)$ the ensemble member:
\begin{equation}\label{upEn}
\theta^{(j)}_{n+1}=\theta^{(j)}_{n}+C^{\theta\omega}_n(C_n^{\omega\omega}+\alpha_n \Gamma)^{-1}(y^{(j)}-\omega^{(j)}_n),
\end{equation}
where the empirical covariances $C^{\theta\omega}_n, C_n^{\omega\omega}$ are given by
$$C^{\theta \omega}_{n} = \frac{1}{N_e-1}\sum^{N_e}_{j=1}(\theta^{(j)}_n-\bar{\theta}_n)(f(\theta^{(j)}_n)-\bar{\omega}_n)^T$$
$$C^{\omega \omega}_{n} = \frac{1}{N_e-1}\sum^{N_e}_{j=1}(f(\theta^{(j)}_n)-\bar{\omega}_n)(f(\theta^{(j)}_n)-\bar{\omega}_n)^T.$$
Here $\bar{\theta}_n$ denotes the average of $\theta^{(j)}_n$ and $\bar{\omega}_n$ denotes the average of $f(\theta^{(j)}_n)$.  The regularizing iterative ensemble Kalman smoother is terminated according to the following discrepancy principle
\begin{equation}
\|\Gamma^{-1/2}(y-\bar{\omega}_n)\| \leq \tau \eta,
\end{equation}
where $\tau$ is a constant and the noise level $\eta$ is defined by
$$\eta = \|\Gamma^{-1/2}(y-f(\theta^{\dag}))\|,$$
here $\theta^{\dag}$ denotes the truth properties.
Note that  the update procedure (\ref{upEn}) with a fixed regularization parameter $\alpha_n=1$ is motivated by the application of Kalman methods for solving Bayesian inference problems when the model $f$ is linear, and the underlying prior distribution is Gaussian \cite{Tarantola2005inverse}. For nonlinear models, the same choice of $\alpha_n$ may lead to instabilities, however. The numerical results of the work  \cite{Iglesias2016regularizing} show that such instabilities can be addressed by choosing  the regularization parameter $\alpha_n$  according to the following criteria
\begin{equation}\label{alrule}
\alpha^N_n \|\Gamma^{1/2}(C_n^{\omega\omega}+\alpha^N_n\Gamma)^{-1}(y^{(j)}-\bar{\omega}_n)\|\geq \rho \|\Gamma^{-1/2}(y^{(j)}-\bar{\omega}_n)\|.
\end{equation}
The detail of the regularizing iterative ensemble Kalman smoother is given in Algorithm \ref{alg:IREKS}.

\begin{algorithm}[t]
  \caption{ Regularizing iterative ensemble Kalman smoother \cite{Iglesias2016regularizing}}
  \label{alg:IREKS}
  \begin{algorithmic}[1]
%   \Procedure{RunChain}{$\widetilde{L}_N,\pi(z),q,m$}
  %  \Require
%    The low-fidelity model $u^L=\sum_{\mb{m}\in\Lambda_N} u^L_{\mb{m}} \Phi_{\mb{m}}(Z)$;  the high-fidelity model $u^H$;   and  the order of $N_C$;
 %  \Ensure
%      Ensemble of classifiers on the current batch, $E_n$;
   \State  \textbf{Prior ensemble and perturbed noise.} Let $\rho<1$ and $\tau\geq1/\rho$. Generate
$$\theta^{(j)}_0\sim \pi(\theta), y^{(j)}=y+\xi^{(j)}, \quad \xi^{(j)}\sim N(0,\Gamma), j=1,\cdots, N_e.$$
\noindent Then  for $n=1,\dots, I_{max}$
    \State \textbf{Prediction step:}  Evaluate
\begin{equation}\label{fevaluate}
\omega^{(j)}_n=f(\theta^{(j)}_n),\quad j=1,\cdots, N_e
\end{equation}
and define $\bar{\omega}_n=\frac{1}{N_e}\sum^{N_e}_{j=1}\omega^{(j)}_n.$
    \State \textbf{Discrepancy principle:} If
\begin{equation}\label{dprule}
\|\Gamma^{-1/2}(y-\bar{\omega}_n)\| \leq \tau \eta,
\end{equation}
stop. Output $\bar{\theta}_n=\frac{1}{N_e}\sum^{N_e}_{j=1}\theta^{(j)}_n.$
  \State \textbf{Analysis step:} Define $C^{\theta \omega}_{n}, C^{\omega \omega}_{n}$ by
$$C^{\theta \omega}_{n} = \frac{1}{N_e-1}\sum^{N_e}_{j=1}(\theta^{(j)}_n-\bar{\theta}_n)(\omega^{(j)}_n-\bar{\omega}_n)^T$$
$$C^{\omega \omega}_{n} = \frac{1}{N_e-1}\sum^{N_e}_{j=1}(\omega^{(j)}_n-\bar{\omega}_n)(\omega^{(j)}_n-\bar{\omega}_n)^T.$$
\noindent Update each ensemble member:
\begin{equation}
\theta^{(j)}_{n+1}=\theta^{(j)}_{n}+C^{\theta\omega}_n(C_n^{\omega\omega}+\alpha_n \Gamma)^{-1}(y^{(j)}-\omega^{(j)}_n), \, j=1,\cdots, N_e,
\end{equation}
where $\alpha_n$ is chosen by the following sequence
\begin{equation}
\alpha^{i+1}_n=2^i \alpha_n^0,
\end{equation}
where $\alpha^0_n$ is an initial guess. We then define $\alpha_n=\alpha^N_n$ where $N$ is the first integer such that
\begin{equation}
\alpha^N_n \|\Gamma^{1/2}(C_n^{\omega\omega}+\alpha^N_n\Gamma)^{-1}(y^{(j)}-\bar{\omega}_n)\|\geq \rho \|\Gamma^{-1/2}(y^{(j)}-\bar{\omega}_n)\|.
\end{equation}
%     \EndProcedure
  \end{algorithmic}
\end{algorithm}

It should be noted that the main computational cost of  regularizing iterative ensemble Kalman smoother  per iteration and per ensemble is due to  Eq. (\ref{fevaluate}). The total cost of an $N_e$ size ensemble of EKI is approximately $N_eJ$ forward model evaluations where $J$ is the total number of iterations; when the model is computationally intensive, which is the case for time-depended partial differential equations, the EKI then becomes prohibitive. It is thus natural to construct a surrogate of the forward model before the data are available. In the next section, we will focus on the polynomial chaos (PC) expansions based surrogate, which is widely used in applied mathematics and engineering.

\subsection{PC-based EKI}
 We first assume that the components of the uncertain parameter vector $\theta=(\theta^{1},\cdots,\theta^{d})$ are mutually independent  and  $\theta^{i}$ has marginal probability density $\pi_i(\theta^{i}): \Theta_i \rightarrow \mathbb{R}^{+}$. Then
$\pi(\theta)=\prod^{d}_{i=1}\pi_i(\theta^{i})$
is the joint probability density of the random vector $\theta$ with the support
$\Theta:=\prod^{d}_{i=1}\Theta_i \in \mathbb{R}^{d}$.

The PC expansion is an orthogonal polynomial approximation to model output $f(\theta)$ which has been broadly used in uncertainty quantification in recent decades \cite{ghanem_and_spanos_book,xiu2010book}. Let $\alpha = (\alpha^1, \cdots, \alpha^d)\in \mathbb{N}^d_0$ be a multi-index with $|\alpha|=\alpha^1+\cdots+\alpha^d$, and $N\geq 0$ be an integer. Then the $N$th degree PC expansion $f_N(\theta)$ of function $f(\theta)$ is defined as
 \begin{eqnarray} \label{gpcexpansion}
f_N(\theta)=\sum_{\alpha \in \Lambda_N^d} c_{\alpha} \Psi_{\alpha}(\theta),  \quad \Lambda_N^d=\{\alpha\in \mathbb{N}^{d}_0: |\alpha|\leq N\}
 \end{eqnarray}
where $\{c_{\alpha}\}$ are the unknown expansion coefficients, and the basis functions $\{\Psi_{\alpha}\}$ are orthonormal under the density $\pi$, that is,
\begin{eqnarray*}
(\Psi_{\alpha},\Psi_{\beta})_{\pi}=\int_{\Gamma}\Psi_{\alpha}(\theta)\Psi_{\beta}(\theta)\pi(\theta)d\theta=\delta_{\alpha,\beta}.
\end{eqnarray*}

By placing an order for the orthogonal polynomials, we can rewrite Eq. (\ref{gpcexpansion}) as the following single index version
\begin{eqnarray} \label{gpce}
f_{N}(\theta)=\sum_{\alpha \in \Lambda_N^d} c_{\alpha}\Psi_{\alpha}(\theta) = \sum^M_{m=1} c_m \Psi_m(\theta),
 \end{eqnarray}
where
\begin{eqnarray}\label{tdterms}
M=\mbox{card} (\Lambda_N^d)= {d+N \choose d}.
\end{eqnarray}
Then the main issue in using PC expansion is to efficiently evaluate the unknown coefficients $\{c_m\}$.  In recent years, more and more attention has been devoted to determine the expansion coefficients based on the data information $\mathcal{D} = \{(\theta_i, f(\theta_i))\}_{i=1}^Q$. In the standard discrete least square method (LSM) \cite{Tang2014discrete,Zhou2015weighted}, we  seek to find the PC coefficients by solving the optimization problem
\begin{eqnarray}
\Big\{c_m\Big\}^M_{m=1} = \arg \min_{c_m} \sum^Q_{i=1}\Big[\Big(f(\theta_i)-\sum^M_{m=1}c_m\Psi_m(\theta_i)\Big)\Big]^2.
\end{eqnarray}
This problem can be written algebraically
\begin{eqnarray}\label{lseq}
\mb{c}^{\#} = \arg \min_{\mb{c} \in \R^M} \|\mb{\Psi c}-\mb{b}\|_2^2,
\end{eqnarray}
where $\mb{c} =(c_1,\cdots, c_M)^T$ denotes the vector of PC coefficients,  $\mb{\Psi}\in \R^{Q\times M}$ denotes the Vandermonde matrix with entries $\mb{\Psi}_{ij}=\Psi_j(\theta_i), \quad i=1,\cdots, Q, j=1,\cdots, M$, and $\mb{b}=(f(\theta_1),\cdots,f(\theta_Q))^T \in \R^Q$ is the vector of samples of $f(\theta)$.   In this paper, we use weighted discrete least square method \cite{Narayan2014} to estimate these coefficients.  Let $\mb{W} =\mbox{diag} (w_1,\cdots,w_Q)$ be a diagonal matrix with positive entries $w_i = \frac{M}{\sum^M_{m=1} \Psi_m^2(\theta_i)}$, a weighted formulation can be written as
\begin{eqnarray}\label{Wlseq}
\mb{c}^{\#} = \arg \min_{\mb{c} \in \R^M} \|\mb{\sqrt{\mb{W}}\Psi c}-\sqrt{\mb{W}}\mb{b}\|_2^2.
\end{eqnarray}

It is clear that after obtaining the approximation of $\mb{c}^{\#}$, one has an explicit functional form $\widetilde{f}_N$. We can then replace the forward model $f$ in (\ref{fevaluate}) by its approximation $\widetilde{f}_N$, and obtain the PC-based EKI algorithm.  Notice that the computational cost of generating $N_e$ samples using $\widetilde{f}_N$ in prediction step requires nothing but sampling of the polynomial expression of (\ref{gpce}) with $N_e$ samples of $\theta$. This cost is minimal because it does not require any simulations of the forward model. Therefore, the main computational cost of PC-based EKI is spent on building the PC representation of the forward model $f$. The total cost of PC-based EKI is approximate $Q>M$ forward evaluations. It should be noted that the accuracy and efficiency of PC depend on the degree of the PC basis. The more PC terms one use, the higher the accuracy one may obtain. However,  as in practice, one can only afford PC expansions with small or moderate PC orders  due to the computational complexity \cite{lu2015JCP,Yan+Zhou2018adaptive}. This can obviously introduce a possibly large model error unless the problem is well represented by a low-order PC. If the model error is large, then there might be a dramatic difference between the inversion results  and the true solution, see Section  \ref{sec:tests}.  To balance accuracy and efficiency, it is desirable to construct a multi-fidelity model to reduce the computational cost of EKIs, namely, one combining  a small number of forward model evaluations  and  a much larger number of  low-order PC model evaluations to construct a multi-fidelity surrogate \cite{Peherstorfer2016survey}.

\section{Adaptive multi-fidelity polynomial chaos approach}\label{sec:method}

\subsection{Multi-fidelity PC based on LSM}

\begin{algorithm}[t]
  \caption{Multi-fidelity PC based on LSM}
  \label{alg:MPC}
  \begin{algorithmic}[1]
%   \Procedure{RunChain}{$\widetilde{L}_N,\pi(z),q,m$}
    \Require
    The low-fidelity model $f^L=\sum_{\alpha\in\Lambda_N^d} u^L_{\alpha} \Psi_{\alpha}(\theta)$;  the high-fidelity model $f^H$;   and  the order of $N_C$;
%    \Ensure
%      Ensemble of classifiers on the current batch, $E_n$;
  \State   Choose $Q=2 {N_C +d \choose d}$ sampling points $\{\theta_i\}$ in the parametric space
    \State Calculate the difference between the $f^H(\theta_i)$ and $f^L(\theta_i)$
    \State Compute the correction PC coefficients $u^C_{\alpha}$ using the least square method
    \State Build the multi-fidelity model by combining $u^L_{\alpha}$ and $u^C_{\alpha}$ using Eq. (\ref{multieq})
%     \EndProcedure
  \end{algorithmic}
\end{algorithm}

In this section, we shall give a brief overview of multi-fidelity polynomial chaos based on LSM. Further details can be found in  \cite{Doostan,Ng2012multi,Palar2016multi}.    The main idea of the multi-fidelity PC approach is to correct the low-fidelity simulation model using a correction term $C:$
\begin{equation}
C(\theta) = f^H(\theta)-f^L(\theta)\approx \sum_{\alpha \in \Lambda_{N_C}^d} u^C_{\alpha} \Psi_{\alpha}(\theta),
\end{equation}
where $f^H$ and $f^L$ are high- and low-fidelity model respectively. Here the unknown coefficients of the additive correction terms $u^C_{\alpha} $ can be calculated by the least squares method. By solving the PC expansions of the correction term, a multi-fidelity model can be approximated via
\begin{equation}
f^H(\theta)= f^L(\theta)+C(\theta)\approx\sum_{\alpha\in\Lambda_N^d} u^L_{\alpha} \Psi_{\alpha}(\theta)+\sum_{\alpha \in \Lambda_{N_C}^d} u^C_{\alpha}\Psi_{\alpha}(\theta),
\end{equation}
where $u^L_{\alpha}$ and $u^C_{\alpha} $ are PC coefficients of the low-fidelity and the correction expansions, respectively.

In practical applications, the indices of correction expansion must be a subset of low-fidelity expansion indices. For example, to construct an $N$-th order multi-fidelity expansion, one can use an $N$-th order low-fidelity PC expansion combined with an $N_C$-th order ($N_C\leq N$) correction expansion.  The multi-fidelity PC expansion can then be expressed as
\begin{equation}\label{multieq}
f^M(\theta)=\sum_{\alpha\in \Lambda_N^d} u^L_{\alpha}\Psi_{\alpha}+\sum_{\alpha\in \Lambda_{N_C}^d} u^C_{\alpha} \Psi_{\alpha}=\sum_{\alpha\in \Lambda_{N_C}^d} (u^L_{\alpha}+u^C_{\alpha}) \Psi_{\alpha}+\sum_{\alpha\in \Lambda_N^d  \backslash \Lambda_{N_C}^d } u^L_{\alpha}\Psi_{\alpha},
\end{equation}
In this way, the multi-fidelity PC introduces an efficient PC approach where the lower-order indices of the low-fidelity PC coefficients are corrected through high-fidelity computations.  The details of the multi-fidelity PC based on LSM are shown in Algorithm \ref{alg:MPC}.

\subsection{Adaptive multi-fidelity PC-based EKI}
 As demonstrated in our previous work \cite{Yan+Zhou2018adaptive}, an accurate multi-fidelity PC surrogate can be  adaptively constructed and refined over a sequence of samples close to the concentrated region of the posterior parameter space. This will significantly improve the accuracy  without a dramatic increase in the computational complexity. Based on this idea, we proposed an adaptive approach integrating the multi-fidelity PC surrogate construction and the EKI. The strategy contains the following steps:

 \textbf{Step 1:} Initialization:  Choosing $Q_1=2 {N+d \choose d}$ sampling points from the prior distribution. Then evaluate the forward model $f$ at these points and build the prior-based PC surrogate $f^L$. Set an initial multi-fidelity surrogate $f^M=f^L$. Generate $N_e$ parameter realizations from the prior distribution as the initial ensemble.

\textbf{Step 2:}  At the $n$-th iteration step, we can generate the system outputs for the ensemble realizations with the surrogate $f^M$.  Update the parameter ensemble with the EKI formula, i.e., Algorithm \ref{alg:IREKS}.

\textbf{Step 3:} Compute the ensemble mean $\bar{\theta}_{n+1}=\frac{1}{N_e}\sum^{N_e}_{j=1}\theta^{(j)}_{n+1}$ and
the following relative error
\begin{equation}\label{conerr}
err = \frac{\|f(\bar{\theta}_{n+1})-f^M(\bar{\theta}_{n+1})\|_{\infty}}{\|f(\bar{\theta}_{n+1})\|_{\infty}}.
\end{equation}
When the relative error $err$ is less than the user-given threshold $tol$, we suppose that the surrogate model is accurate enough and thus it is used directly in EKI. If the error indicator $err$ exceeds   $tol$, we shall refine the multi-fidelity model $f^M$ using Algorithm \ref{alg:MPC}. In particular, we shall choose $Q_2=2 {N_C+d \choose d}$ random points $\{z^{(i)}\}$ in a ball centered at $\bar{\theta}_{n+1}$,  i.e., $z^{(i)}\in B(\bar{\theta}_{n+1},R):=\Big\{z: \|z-\bar{\theta}_{n+1}\|_{\infty} \leq R\Big\}$, to perform the true model evaluations and then construct a new multi-fidelity model via (\ref{multieq}).

\textbf{Step 4:} Repeat Steps 2-3 until one of the stop criteria of EKI is met.

For the present application, the cost of Step 2  is negligible compared to the cost of the update the multi-fidelity PC surrogate, i.e. Step 3. The total cost of adaptive multi-fidelity PC (AMPC) based EKI algorithm is around $(J_1+Q_1+J_2Q_2)$ forward model evaluations where $J_1$ is the number of iterations to converge, $J_2\leq J_1$ is the number of adaptively.  For the forward models considered in Section \ref{sec:tests}, our numerical results  indicate that $J_2$ is typically between 2 and 6 iterations. Thus, for high-dimensionally parameter models, the computational efficiency of the new algorithm may be comparable to  the PC-based EKI  with a large order PC basis. Usually, the total number of the forward model evaluation of AMPC-based EKI is also much smaller than the standard EKI with a large ensemble (e.g. $10^2\sim10^3$). Thus, the computational cost of AMPC-based EKI is also significantly reduced compared to that of the standard EKI. This is will be demonstrated in the following numerical experiments.

\section{Numerical Examples}\label{sec:tests}

In this section, we present a two dimensional time fractional  PDE inverse problem to illustrate the accuracy and efficiency of the adaptive multi-fidelity PC approach. To better present the results, we shall perform the following three-types of approaches:
\begin{itemize}
\item The {\it conventional EKI}, or the direct EKI approach based on the forward model evaluations.
\item The EKI approach based on a prior-PC surrogate model evaluations.
\item The AMPC approach presented in Section \ref{sec:method}.
\end{itemize}
In our figures and results, we will use ``Direct" to denoted the conventional EKI, ``PC" to denoted the PC-EKI, and ``AMPC" to denote the AMPC algorithm.  Since EKI is a Monte Carlo-based method, its performance is affected by the specific initial ensemble, especially when the ensemble size is small.  In order to illustrate the effect of the initial ensemble, we will show the output for 50 different initializations and report the mean results along with $20\%$ and $80\%$ quantiles, for each example.  We will also plot the final iteration reconstruction arising from one of those initializations.  All the computations are performed using MATLAB 2015a on an Intel-i5 desktop computer.

\subsection{Problem setup}
Consider the following two dimensional time-fractional PDEs in the physical domain  $\Omega =[0, 1]^2$
\begin{eqnarray}\label{2dtfpde}
\begin{array}{rl}
^cD_t^{\alpha}u-\nabla\cdot(\kappa(x;\theta) \nabla u(x,t))&=e^{-t}\exp\Big(-\frac{\|x-(0.25,0.75)\|^2}{2\times0.1^2}\Big),\quad \Omega\times [0, 1],\\
\nabla u \cdot \textbf{n}&=0, \quad \mbox {on} \,\partial{\Omega},\\
u(x,0)&=0, \quad  \mbox{in}\, \Omega.
 \end{array}
\end{eqnarray}

The goal is to determine the permeability field $\kappa(x;\theta)$ from noisy measurements of the $u$-field at a finite set of locations and times.  Here $^cD_t^{\alpha} (0<\alpha<1)$ denotes the Caputo fractional derivative of order $\alpha$. In the numerical simulation, we  solve the equation (\ref{2dtfpde}) using a finite difference/ spectral approximations (\cite{Lin+Xu2007}) with time step $\Delta t=0.01$ and polynomial degree $P=6$.  In order not to commit an 'inverse crime', we generate the data by solving the forward problem using a higher order (P=10) than that is used in the inversion.

In order to measure the accuracy of the numerical approximation $\bar{\kappa}$ with respect to the exact solution $\kappa^{\dag}$, we use the relative error $rel(\kappa)$ defined as
\begin{eqnarray*}
rel(\kappa)=\frac{\sqrt{\sum^{N_0}_{i=1}(\bar{\kappa}_i-\kappa^{\dag}_i)^2}}{\sqrt{\sum^{N_0}_{i=1}(\kappa^{\dag}_i)^2}},
\end{eqnarray*}
where $\bar{\kappa}_i$ and $\kappa^{\dag}_i$ are the numerical  and exact solutions evaluated at the $i$th node, respectively. Here $N_0$ is the total number of resolution points.  In all our tests, unless otherwise specified, we shall use the following parameters  $\alpha = 0.5,  \rho=1/\tau=0.7, I_{max}=30, R=0.2$.

\subsection{Example 1: a nine-dimensional inverse problem}
\begin{figure}
\begin{center}
  \begin{overpic}[width=.45\textwidth,trim=20 0 20 15, clip=true,tics=10]{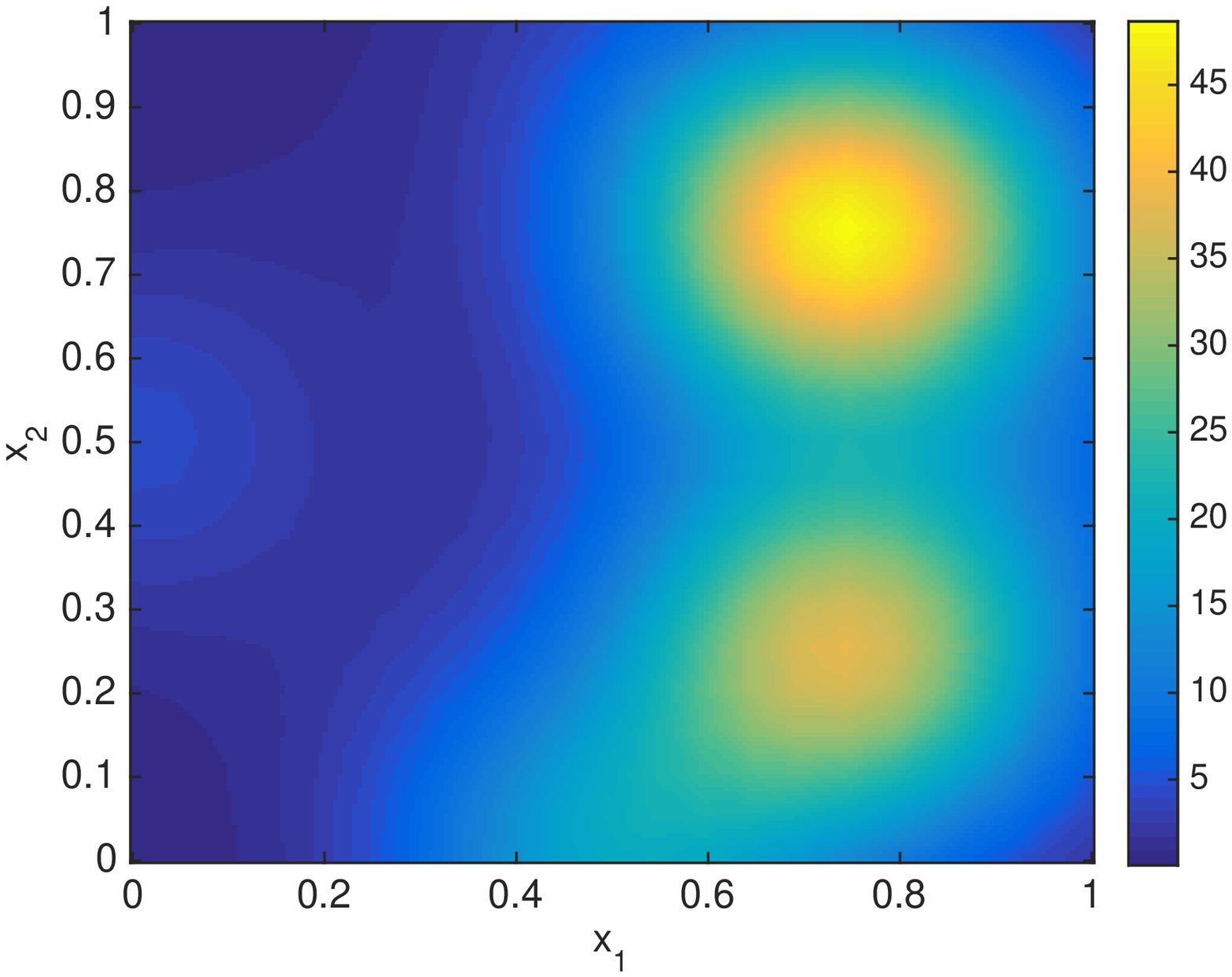}
  \end{overpic}
   \begin{overpic}[width=.45\textwidth,trim=20 0 20 15, clip=true,tics=10]{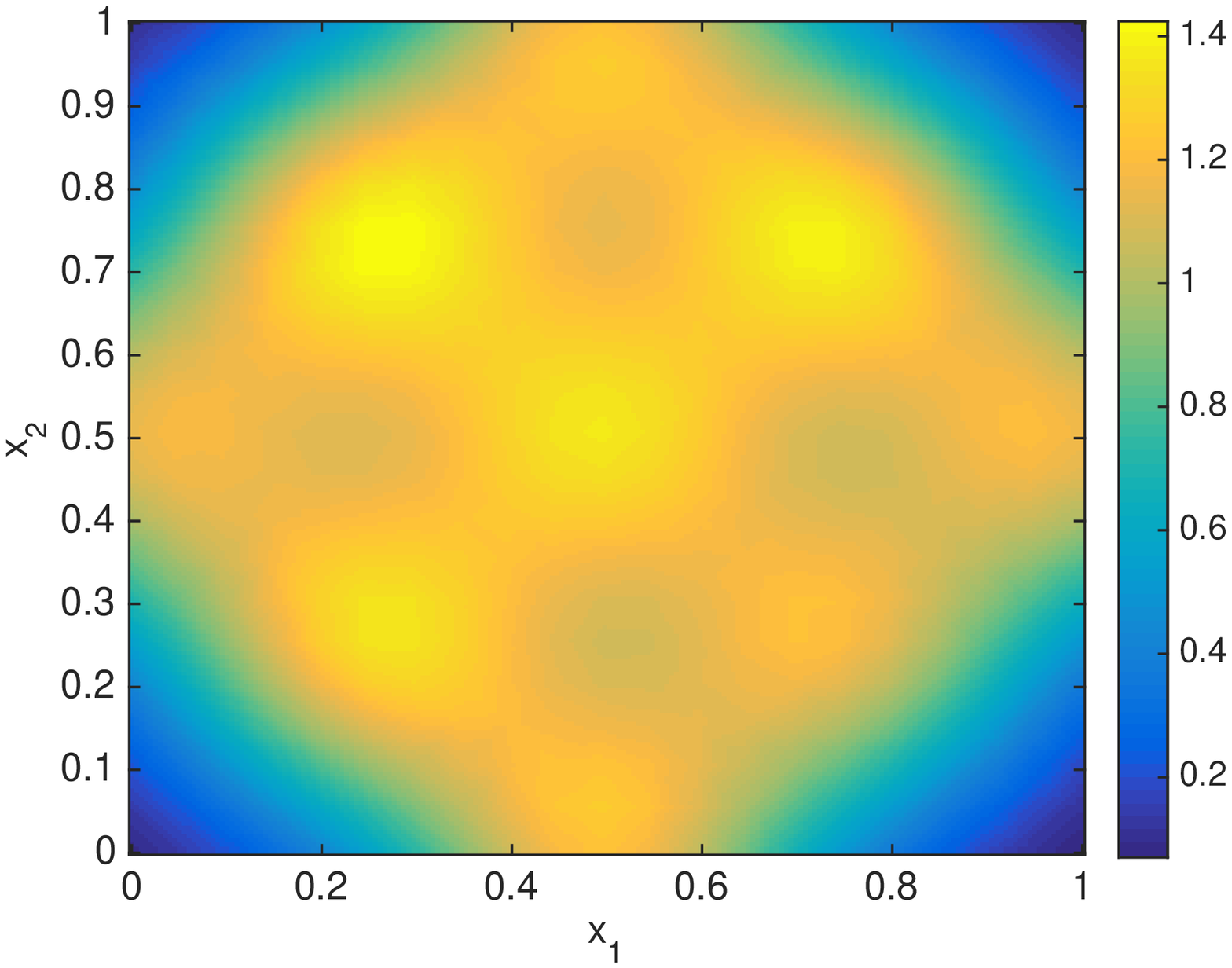}
  \end{overpic}
\end{center}
\caption{Example 1. Left: the true permeability used for generating the synthetic data sets. Right: the initial ensemble mean.}\label{exact_eg1}
\end{figure}

In this example, we consider the following permeability field $\kappa(x;\theta)$
\begin{eqnarray*}
\kappa(x; \theta)=\sum^{9}_{i=1}\theta^i \exp(- 0.5\frac{\|x-x_{0,i}\|^2}{0.15^2}),
\end{eqnarray*}
where $\{x_{0,i}\}^{9}_{i=1}$ are the centers of the radial basis function.  The prior distributions on each of the weights $\theta^i, i=1,\cdots, 9$ are independent and log-normal; that is, $\log(\theta^i)\sim N(0,1)$.   The true permeability field used to generate the test data is shown in Fig.\ref{exact_eg1}. In this example, the true parameter is drawn from  $\log(\theta^i)\sim U(-4,4)$.  The simulation  data are generated by selecting  the values of the states at a uniform $5\times 5$ sensor network. At each sensor location, three measurements are taken at time $t=\{0.25,0.75,1\}$, which corresponds to a total of 75 measurements.  For simplicity, we consider a diagonal measurement error covariance $\Gamma=\sigma^2 I$.  Then the synthetic data $y$ is generated by
\begin{eqnarray*}
y_j=u(x_j, t_j; \theta^{\dag})+\xi_j,
\end{eqnarray*}
with $\xi_j\sim N(0, \sigma^2)$.  In this case, the parameters are far from the prior, and one cannot guarantee the accuracy of prior-based PC approach due to the lack of global accuracy of the PC surrogate.

 \begin{figure}
\begin{center}
  \begin{overpic}[width=0.3\textwidth,trim=20 0 20 15, clip=true,tics=10]{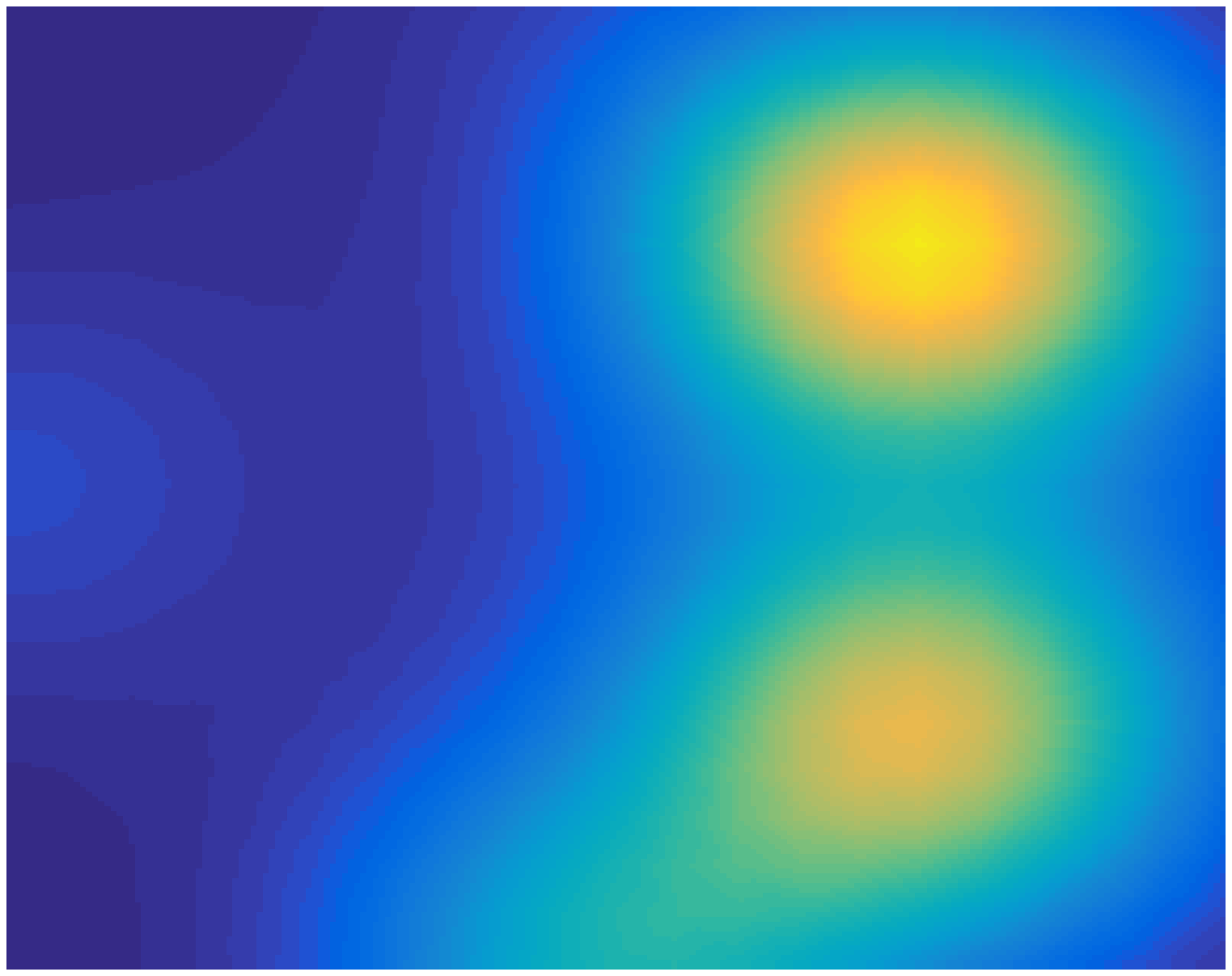}
  \end{overpic}
    \begin{overpic}[width=0.3\textwidth,trim= 20 0 20 15, clip=true,tics=10]{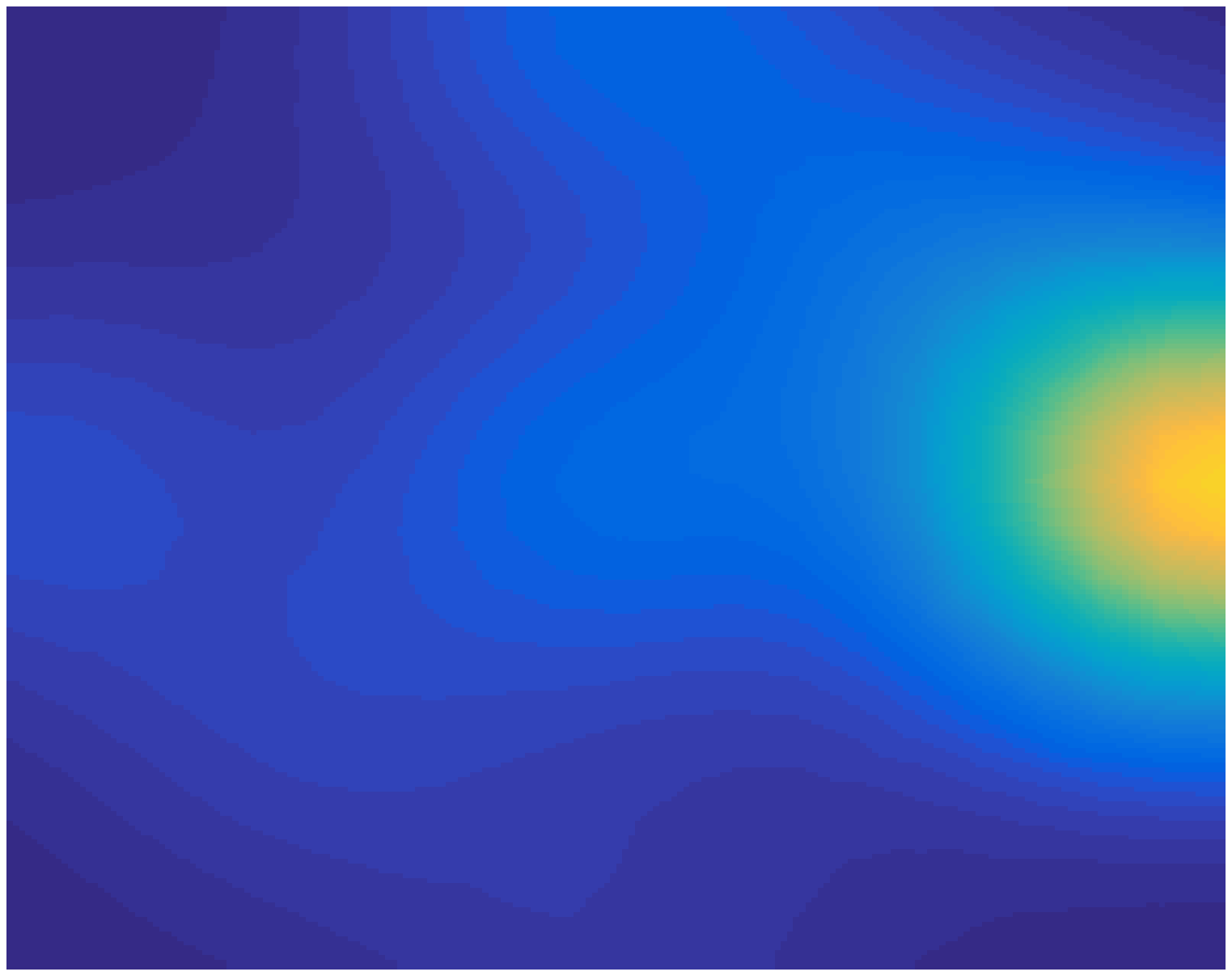}
  \end{overpic}
    \begin{overpic}[width=0.3\textwidth,trim= 20 0 20 15, clip=true,tics=10]{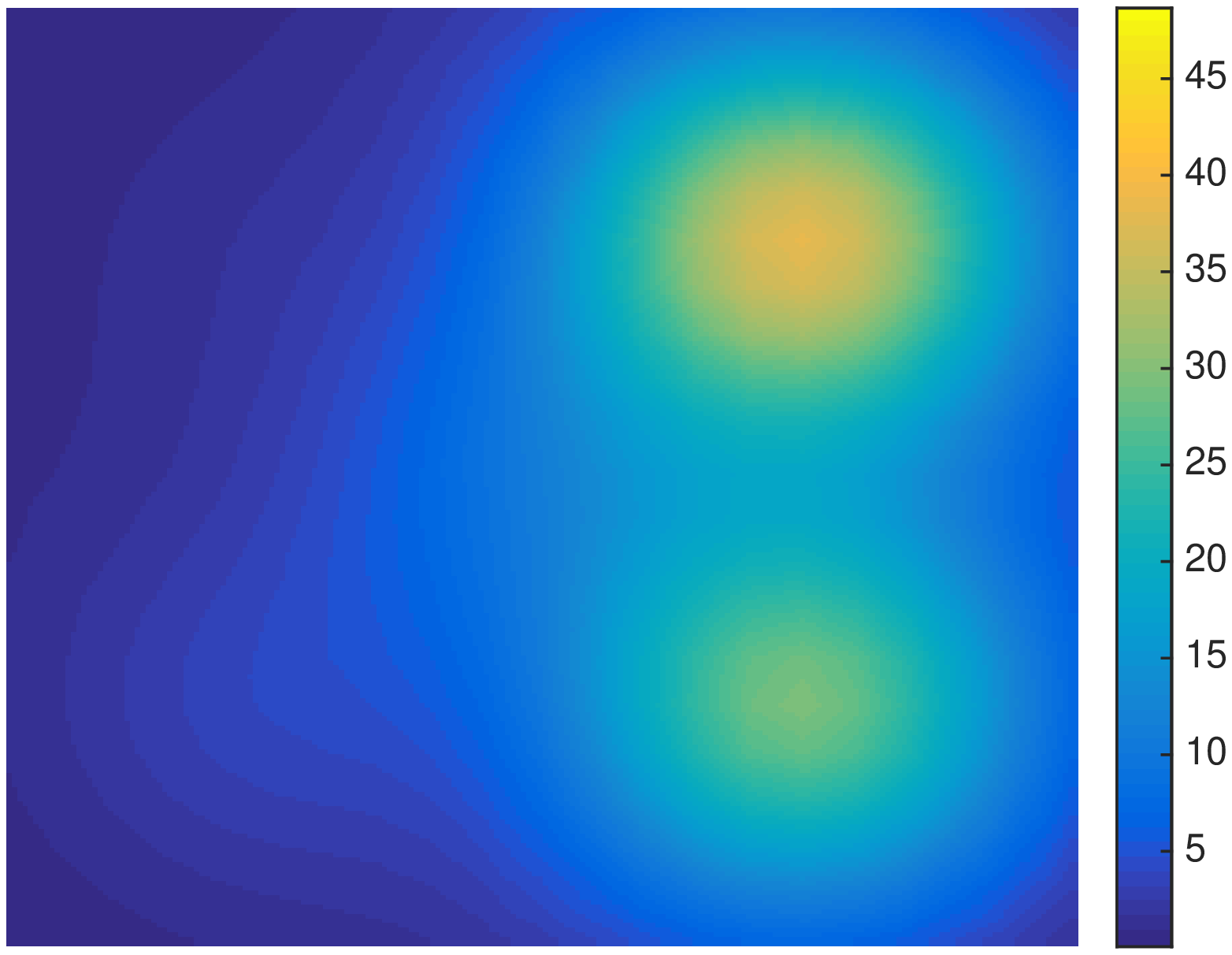}
  \end{overpic}
    \begin{overpic}[width=0.3\textwidth,trim=20 0 20 15, clip=true,tics=10]{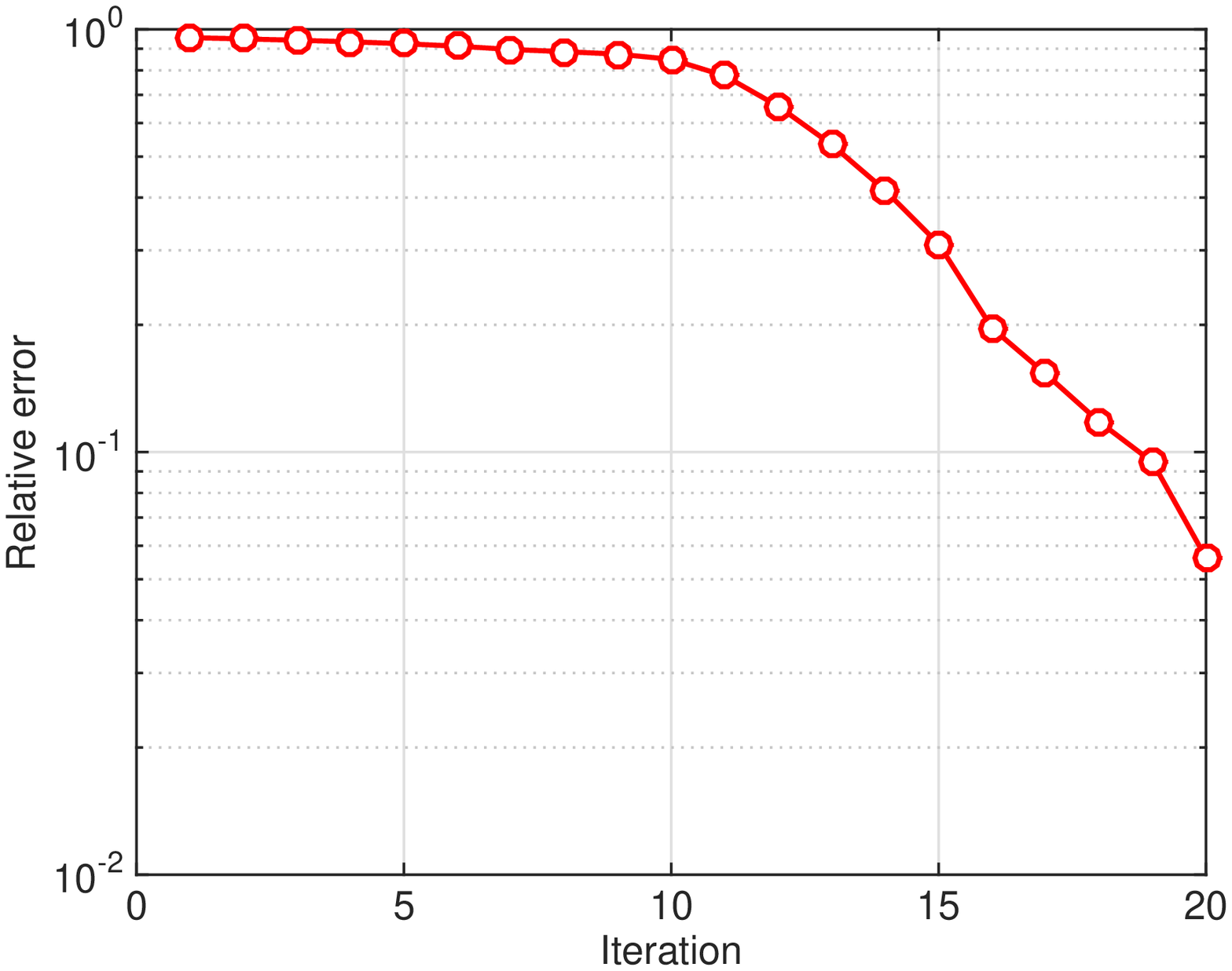}
  \end{overpic}
    \begin{overpic}[width=0.3\textwidth,trim= 20 0 20 15, clip=true,tics=10]{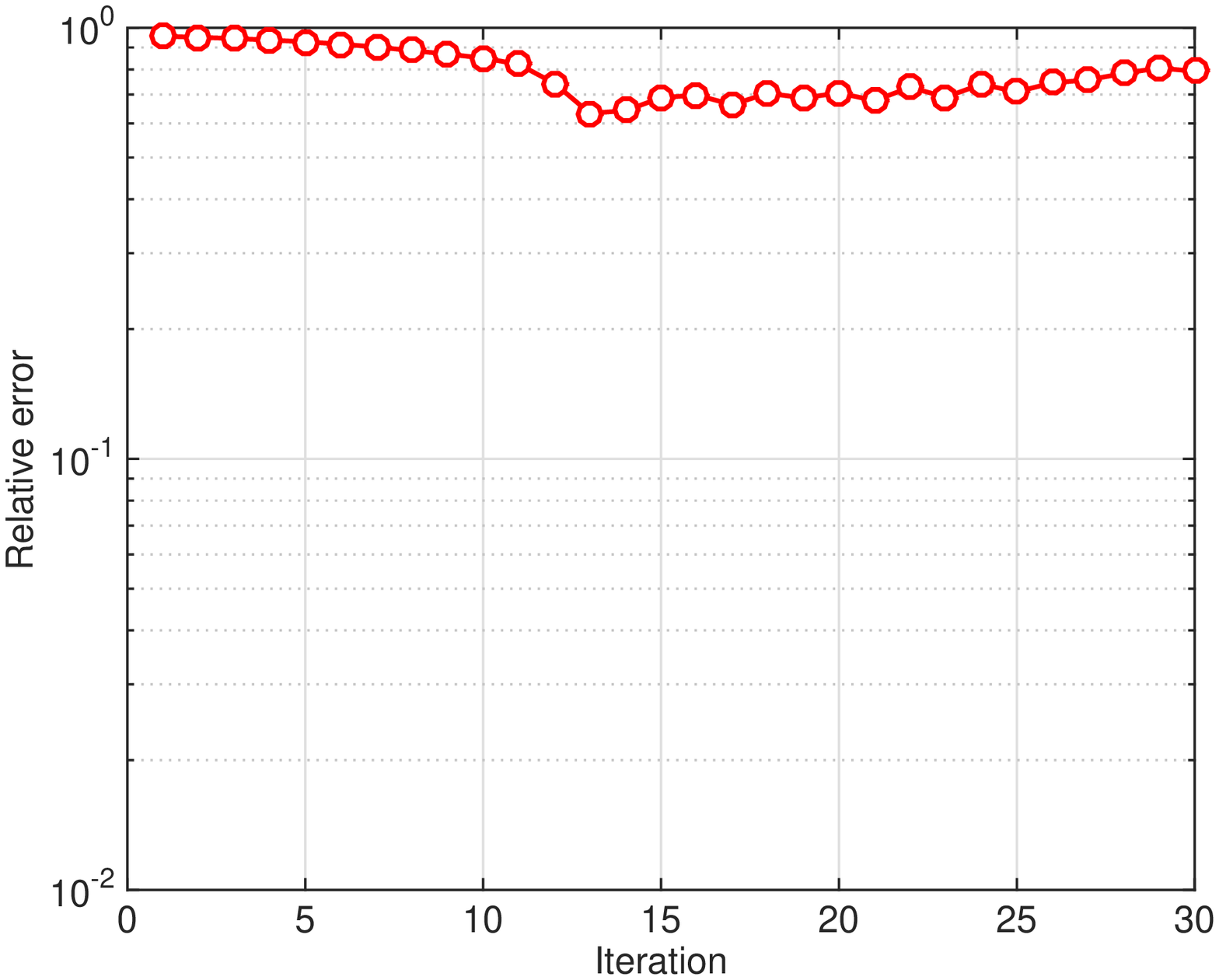}
  \end{overpic}
    \begin{overpic}[width=0.3\textwidth,trim= 20 0 20 15, clip=true,tics=10]{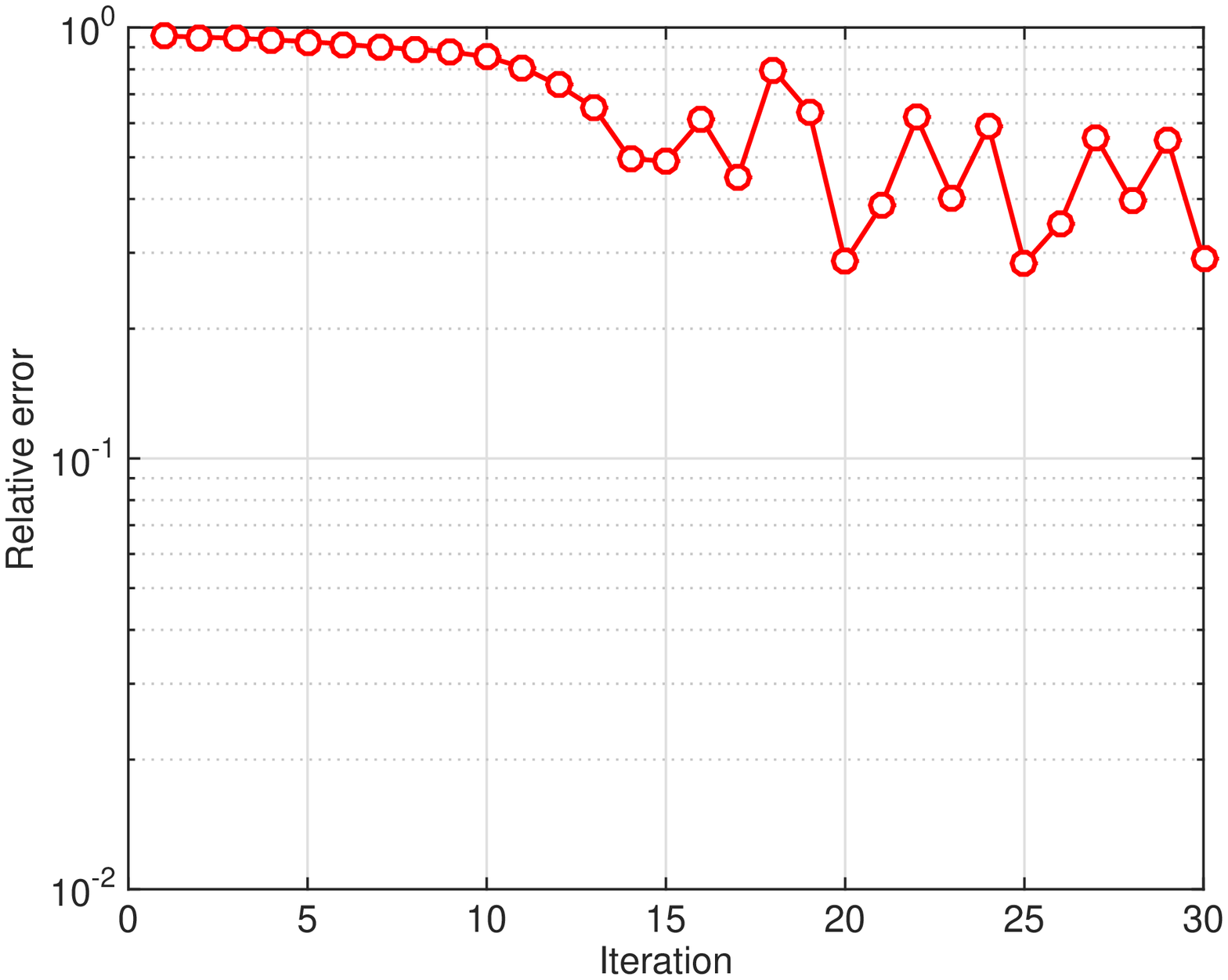}
  \end{overpic}
\end{center}
\caption{Example 1:  Numerical results for the final iteration using $N_e =100$: (Left) Direct; (Middle) PC (N=4); (Right) PC (N=6).}\label{mean-eg1-PC}
\end{figure}

 \begin{figure}
\begin{center}
  \begin{overpic}[width=0.32\textwidth,trim=20 0 20 15, clip=true,tics=10]{figure/direct_9d_eg2.eps}
  \end{overpic}
%    \begin{overpic}[width=0.24\textwidth,trim= 20 0 20 15, clip=true,tics=10]{figure/PC_9d_eg2_N6.eps}
%  \end{overpic}
    \begin{overpic}[width=0.32\textwidth,trim= 20 0 20 15, clip=true,tics=10]{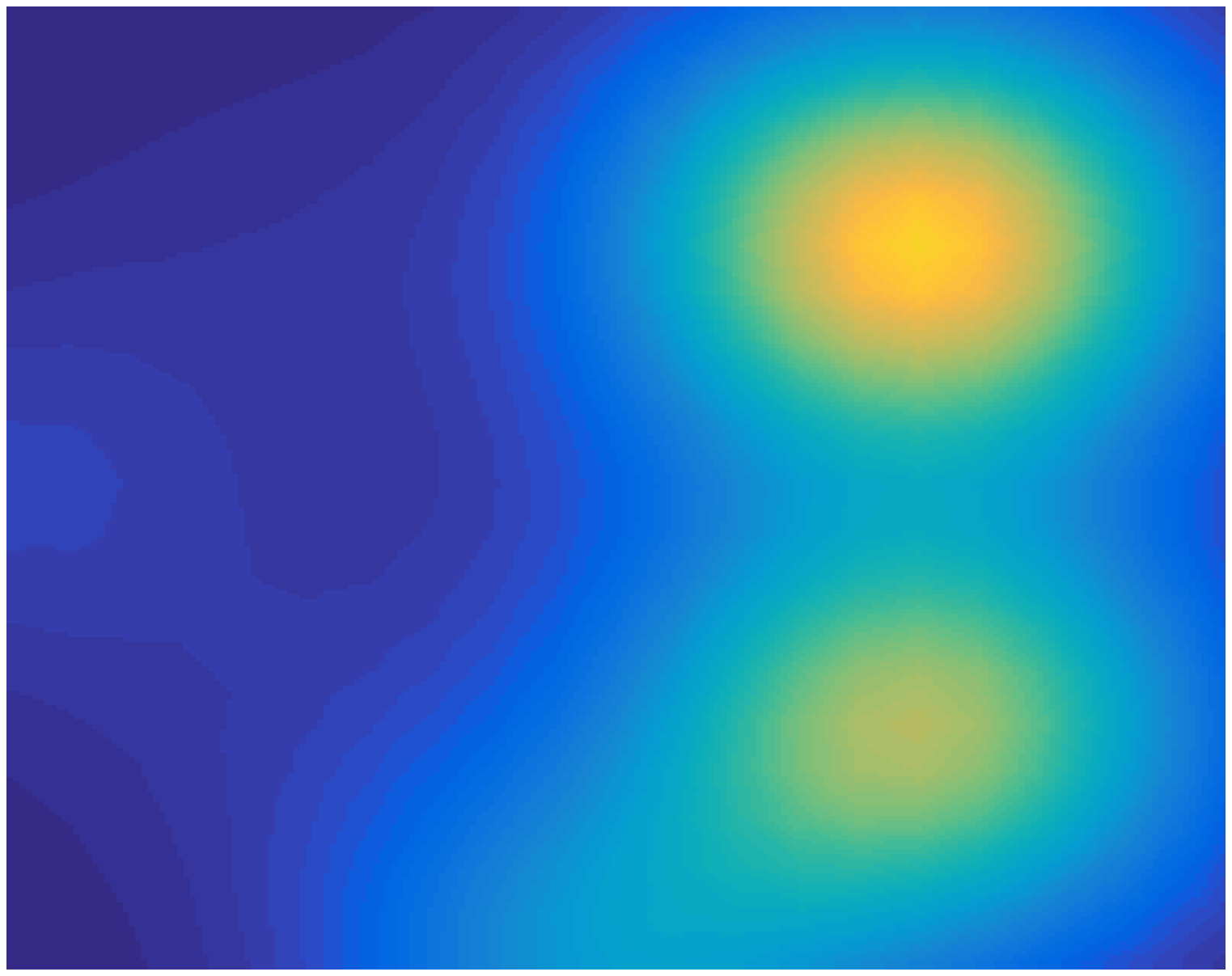}
  \end{overpic}
    \begin{overpic}[width=0.32\textwidth,trim= 20 0 20 15, clip=true,tics=10]{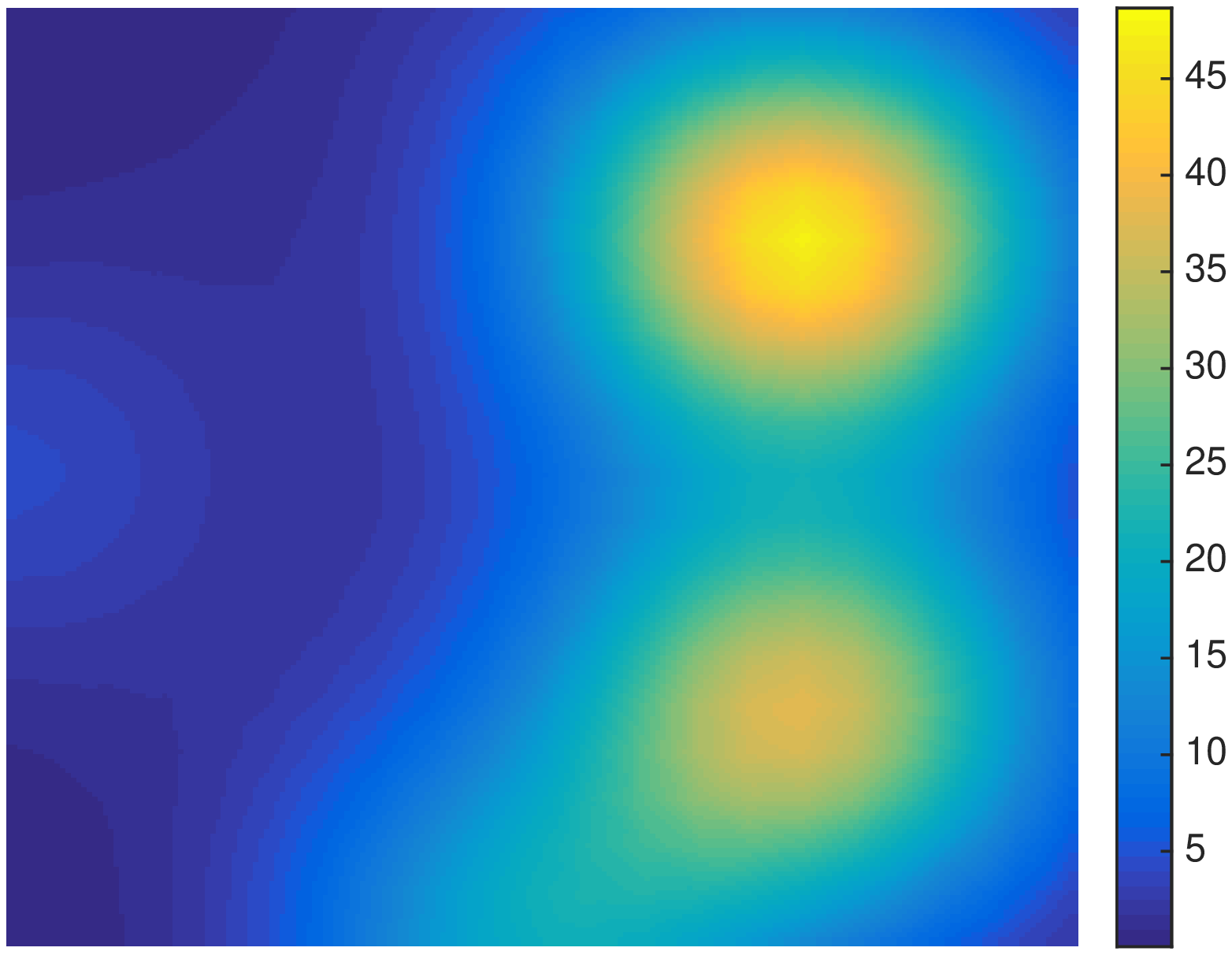}
  \end{overpic}
    \begin{overpic}[width=0.32\textwidth,trim=20 0 20 15, clip=true,tics=10]{figure/direct_9d_eg2_iter.eps}
  \end{overpic}
%    \begin{overpic}[width=0.24\textwidth,trim= 20 0 20 15, clip=true,tics=10]{figure/PC_9d_eg2_N6_iter.eps}
%  \end{overpic}
    \begin{overpic}[width=0.32\textwidth,trim= 20 0 20 15, clip=true,tics=10]{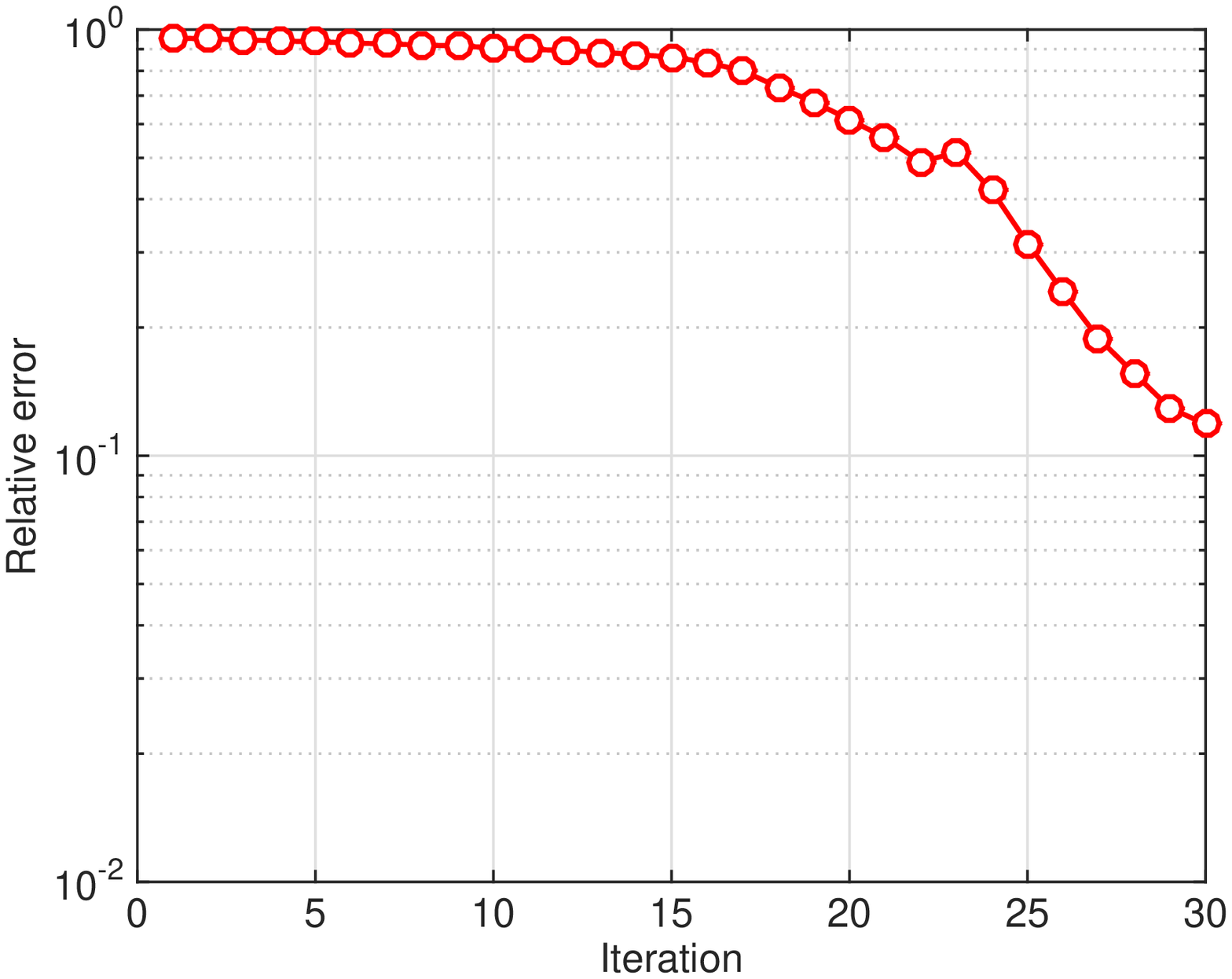}
  \end{overpic}
    \begin{overpic}[width=0.32\textwidth,trim= 20 0 20 15, clip=true,tics=10]{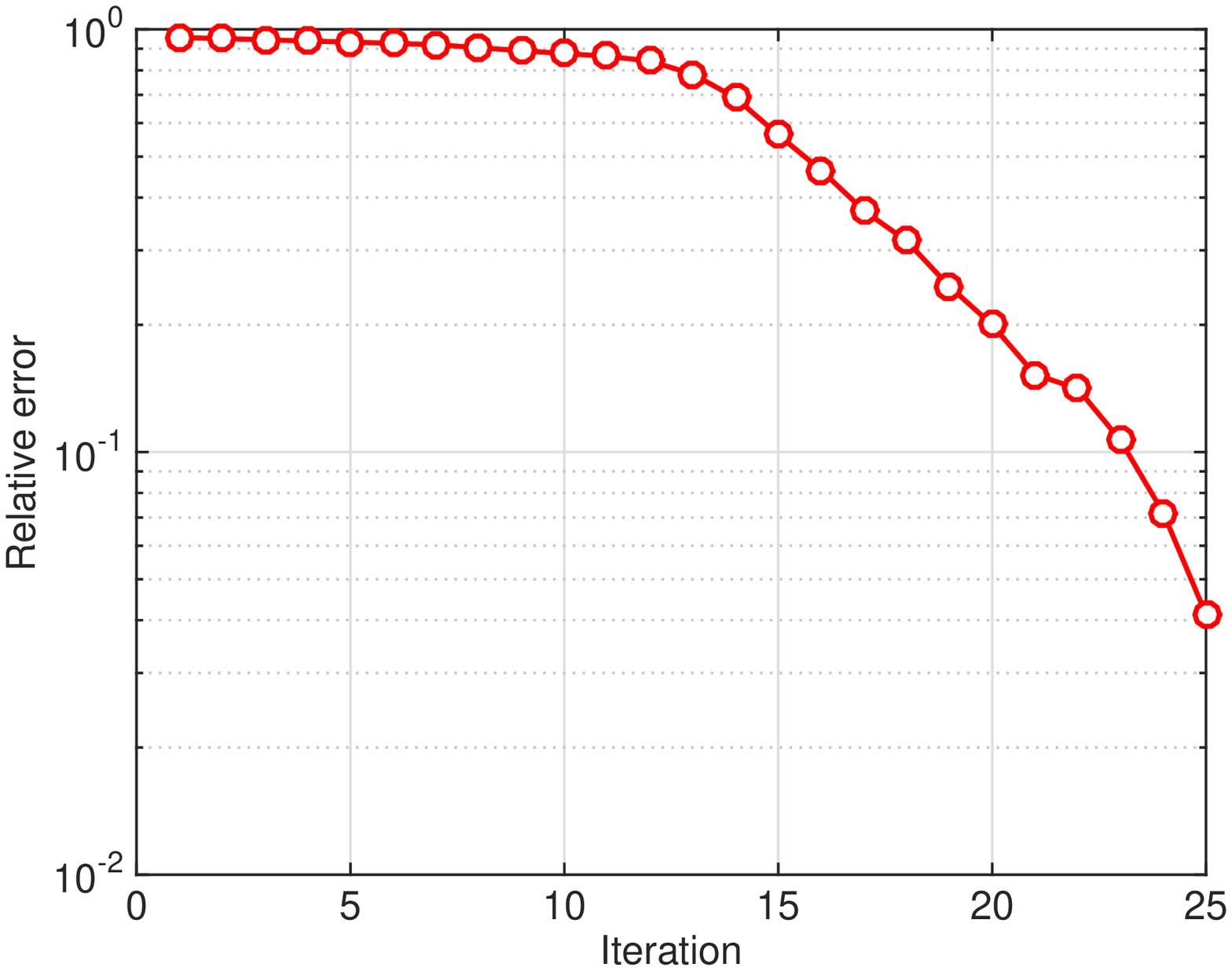}
  \end{overpic}
\end{center}
\caption{Example 1:  Numerical results for the final iteration using $N_e =100$: (Left) Direct; (Middle)  AMPC (N=2, $tol=1\times 10^{-2}$); (Right) AMPC (N=2, $tol=1\times 10^{-3}$).}
\label{mean-eg1-AMPC}
\end{figure}

 \begin{table}[tp]
      \caption{Example 1. Computational times, in seconds, given by three different methods. $N= N_C=2, N_e=100$. }\label{eg1_time}
  \centering
  \fontsize{6}{12}\selectfont
  \begin{threeparttable}
    \begin{tabular}{ c cccccc}
  \toprule
 & \multicolumn{2}{c}{Offline}&\multicolumn{2}{c}{Online}\cr
\cmidrule(lr){2-3} \cmidrule(lr){4-5}

  \multirow{1}{*}{Method}  &$\text{$\#$ of model eval.}$&CPU(s) &$\text{$\#$ of model eval.}$&CPU(s)     &\multirow{1}{*}{Total time(s)}&\multirow{1}{*}{rel}\cr
  \midrule
    Direct                                       & $-$       & $-$          & 2000       &56.71                   & 56.71  & 0.0461\cr
   PC, $N=6$                               & 10010  & 336.59     & $-$          &3.15                     & 339.74  &0.2892\cr
   PC, $N=4$                               & 1430   & 40.25       & $-$          &0.82                      & 41.07  & 0.7921\cr
   AMPC, tol=1e-2    & 110     & 3.69            & 250         &6.28                 & 9.97  &0.1186\cr
   AMPC, tol=1e-3    & 110     & 3.69            & 575         &10.92                 & 14.61  &0.0382\cr
      \bottomrule
      \end{tabular}
    \end{threeparttable}

\end{table}

\begin{figure}
\begin{center}
  \begin{overpic}[width=.45\textwidth,trim=20 0 20 15, clip=true,tics=10]{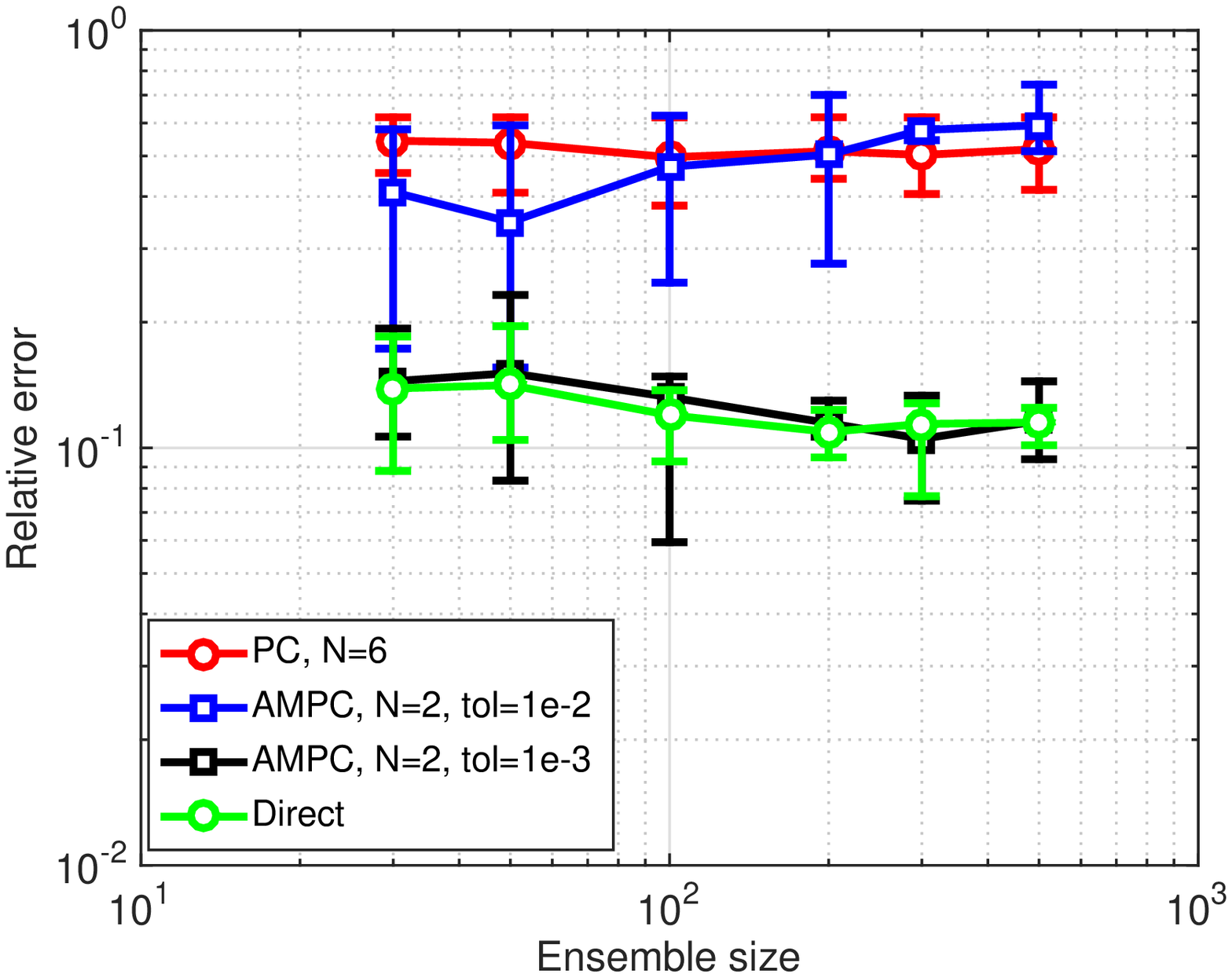}
  \end{overpic}
    \begin{overpic}[width=.45\textwidth,trim=20 0 20 15, clip=true,tics=10]{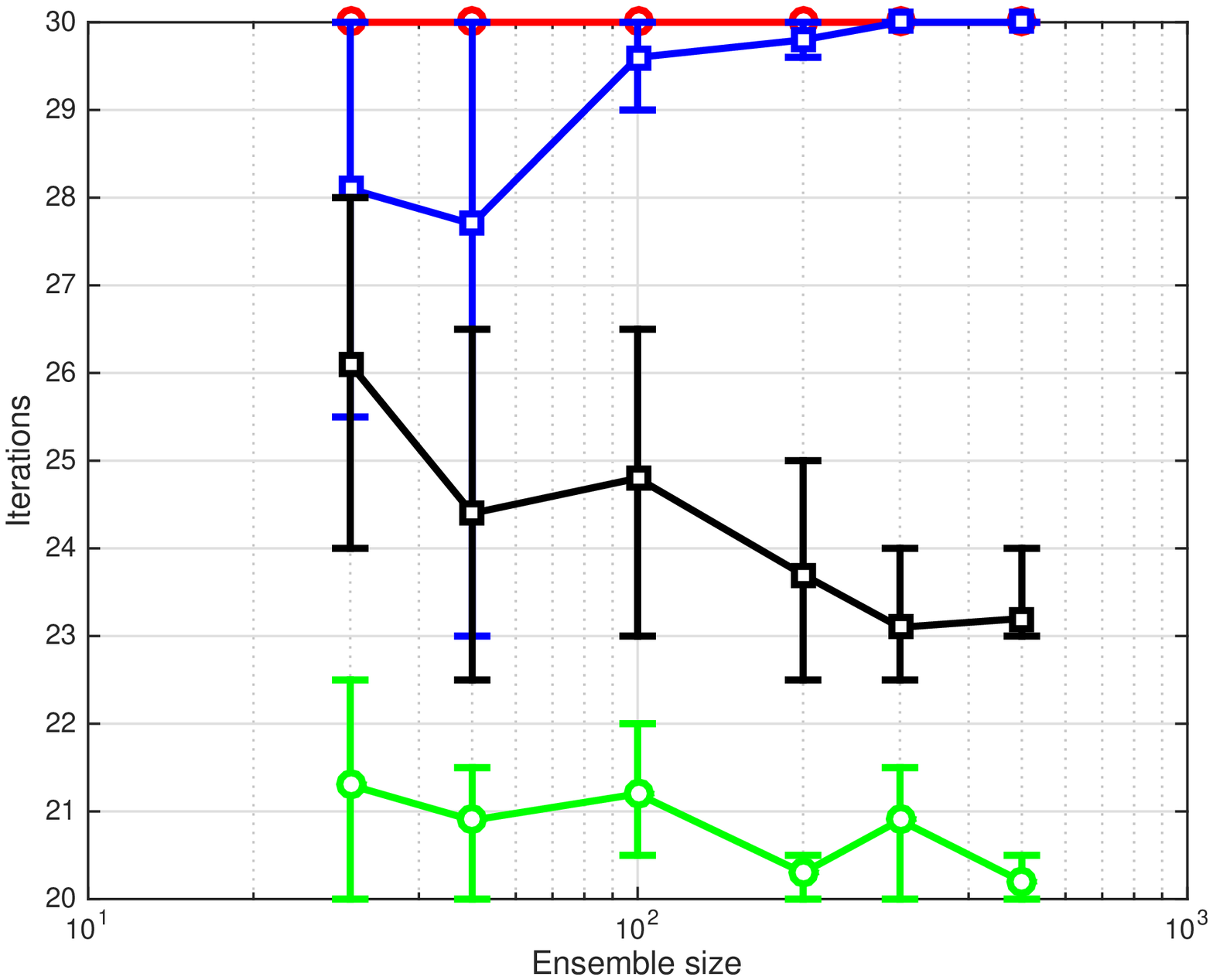}
  \end{overpic}
\end{center}
\caption{Example 1. Numerical results with different ensemble size. Left: values of $rel$ for the estimated results.  Right: number of required iterations.}\label{eror_eg1}
\end{figure}

\begin{figure}
\begin{center}
   \begin{overpic}[width=.45\textwidth,trim=20 0 20 15, clip=true,tics=10]{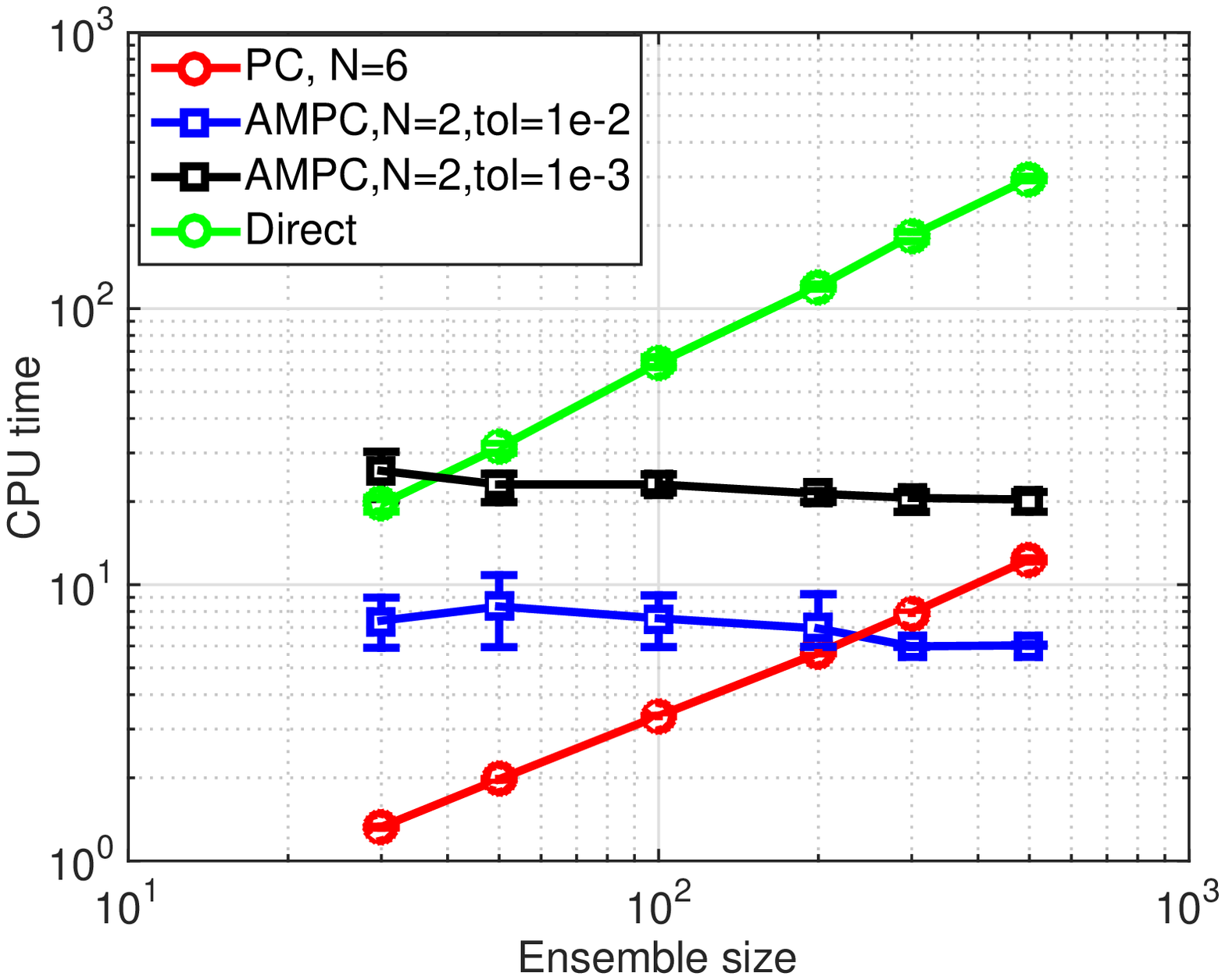}
  \end{overpic}
  \begin{overpic}[width=.45\textwidth,trim=20 0 20 15, clip=true,tics=10]{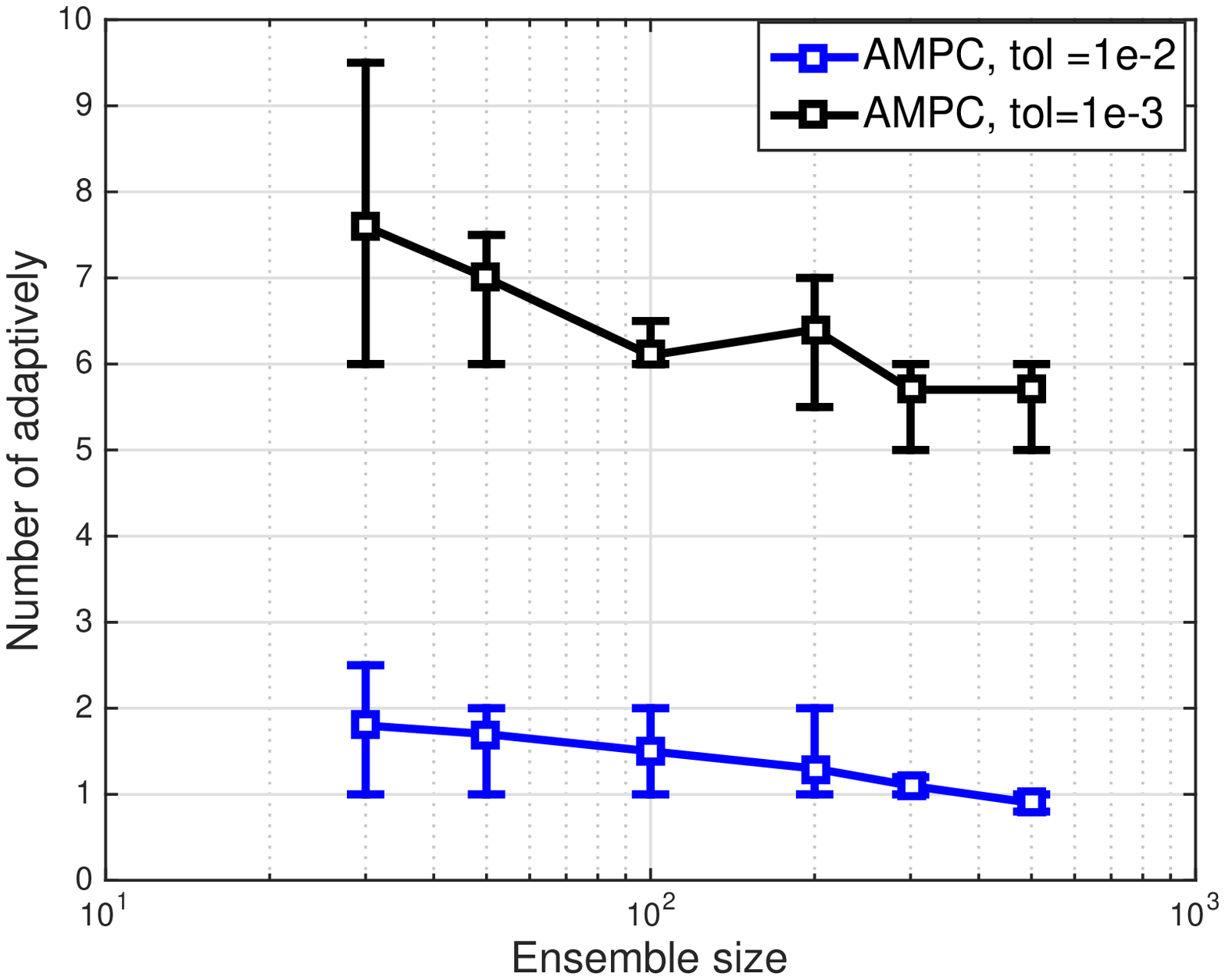}
  \end{overpic}
\end{center}
\caption{Example 1. Numerical results with different ensemble size. Left: online CPU times. Right: number of required adaptively.}\label{cpu_eg1}
\end{figure}

We first investigate the performance of the AMPC method. For our EKI methods, we choose $N_e=100$ ensemble members with the standard deviation $\sigma =1\times 10^{-3}$. The initial mean ensemble is shown in Figure \ref{exact_eg1}. In Figure \ref{mean-eg1-PC}, we show the final iteration reconstruction arising from the conventional EKI,  and the PC-based  EKI. The relative errors $rel$ as a function of the number of iteration also shown in Figure \ref{mean-eg1-PC}.  Since the exact parameter is far from what is assumed in the prior, it is evident from the figures that the results using the PC-based EKI  give a large error.  The relative error of $rel$ was found to be  $0.7921$ and $0.2892$ for the order of PC $N=4$ and $N=6$, respectively, both  much larger than the results by the conventional EKI ($0.0461$). The CPU time of evaluating the conventional EKI is about  56.71s, while the CPU time of PC-based EKI with $N=4$ is about 0.82s. Thus, although using a  PC surrogate in the EKI can gain computational efficiency,  the estimation accuracy cannot be guaranteed.  To improve this one can increase the PC order $N$. However, when the order increases, the cost of constructing the PC surrogate becomes increasingly expensive. For example, constructing a PC model with $N=6$ requires $10,010$ offline model evaluations, and need an offline CPU time of 336.59s.

The corresponding results using AMPC are shown in Figure \ref{mean-eg1-AMPC}. It is not surprising that even a lower PC order $N=2$ is used for AMPC, a rather accurate result can be obtained. As shown in Figure \ref{mean-eg1-AMPC}, the final iteration reconstruction by using the conventional EKI and the AMPC algorithm are almost identical, which indicates the accuracy of the AMPC algorithm.

The computational costs and the relative errors $rel$ of the final iteration, given by three different  algorithms are shown in Table \ref{eg1_time}. The main computational time in the PC-based EKI is the offline model evaluations. The number of such high-fidelity model evaluations with $N=\{6, 4\}$ are $10010$ and 1430, respectively.  Upon obtaining the PC surrogate, the online simulation is very cheap as it does not require any forward model evaluations. For the AMPC, we do need the online high-fidelity model simulations to refine the multi-fidelity PC surrogate.  Nevertheless, in contrast to $2000$  model evaluations in the conventional EKI, the number of model evaluations for the AMPC with  $N=2, tol=\{10^{-2},10^{-3}\}$ are $250$ and $575$, respectively. As can be seen from the last two columns of Table \ref{eg1_time}, the AMPC approach can improve significantly the accuracy, yet without a dramatic increase in the computational time compared to the PC-based EKI.  This demonstrated that the AMPC is more efficient than the PC-based EKI for solving problems which the data contain information beyond what is assumed in the prior.

Next, we consider the effect of the number of ensemble size on the performance of the algorithms. On the left of Figure \ref{eror_eg1}, the $rel$ values for estimated results as a function of ensemble size are plotted. From this figure, we can find that the $rel$ values decrease to a small and stable level as the ensemble size increases for conventional EKI and AMPC with small $tol=1\times 10^{-3}$. Furthermore, the mean number of required iterations also decreases as the ensemble size increases for conventional EKI and AMPC with small $tol=1\times 10^{-3}$, see  the right of Figure \ref{eror_eg1}.  This is because a larger ensemble size guarantees more accurate sensitivity information, resulting in an accuracy numerical results. However, the computational costs will increase as the ensemble size increases. On the left of Figure \ref{cpu_eg1}, the online CPU times of three different algorithms are plotted against increasing numbers of ensemble size.  Interestingly, the online CPU times of the AMPC are almost unchanged as the ensemble size increases.  The reason is  that the main computational cost of AMPC spent  on  refining the multi-fidelity PC model.  From the right of  Figure \ref{cpu_eg1}, we can see that the number of the adaptively for AMPC is almost unchanged. It is about 2 ($tol=1\times 10^{-2}$) or 7 ($tol=1\times 10^{-3}$) iterations.  On the other hand, the prediction steps of EKI are calculated from a large number of realizations generated by the lower order (e.g., $N=2$) multi-fidelity PC model with virtually no additional computational cost. It should be noted that, when the PC order is larger (e.g., N=6), the online CPU times of PC-based EKI will increase when the ensemble size increases.

\subsection{Example 2}
\begin{figure}
\begin{center}
  \begin{overpic}[width=.45\textwidth,trim=20 0 20 15, clip=true,tics=10]{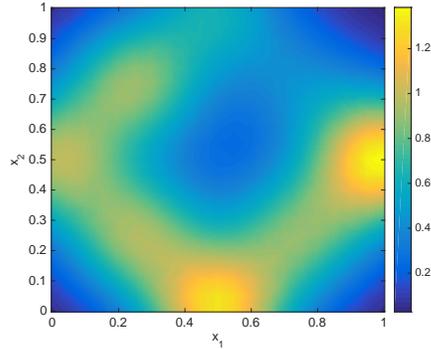}
  \end{overpic}
\end{center}
\caption{Example 2. The true permeability used for generating the synthetic data sets. }\label{exact_eg2}
\end{figure}

As the second example, we define the exact permeability denoted by $\kappa^{\dag}$ and displayed in Fig. \ref{exact_eg2}.  In this example, the true parameter $\theta$ is a draw from the prior distribution described in Example 1. In other words, we consider the best-case-scenario where our prior knowledge includes  the truth.

\begin{figure}
\begin{center}
  \begin{overpic}[width=0.24\textwidth,trim=20 0 20 15, clip=true,tics=10]{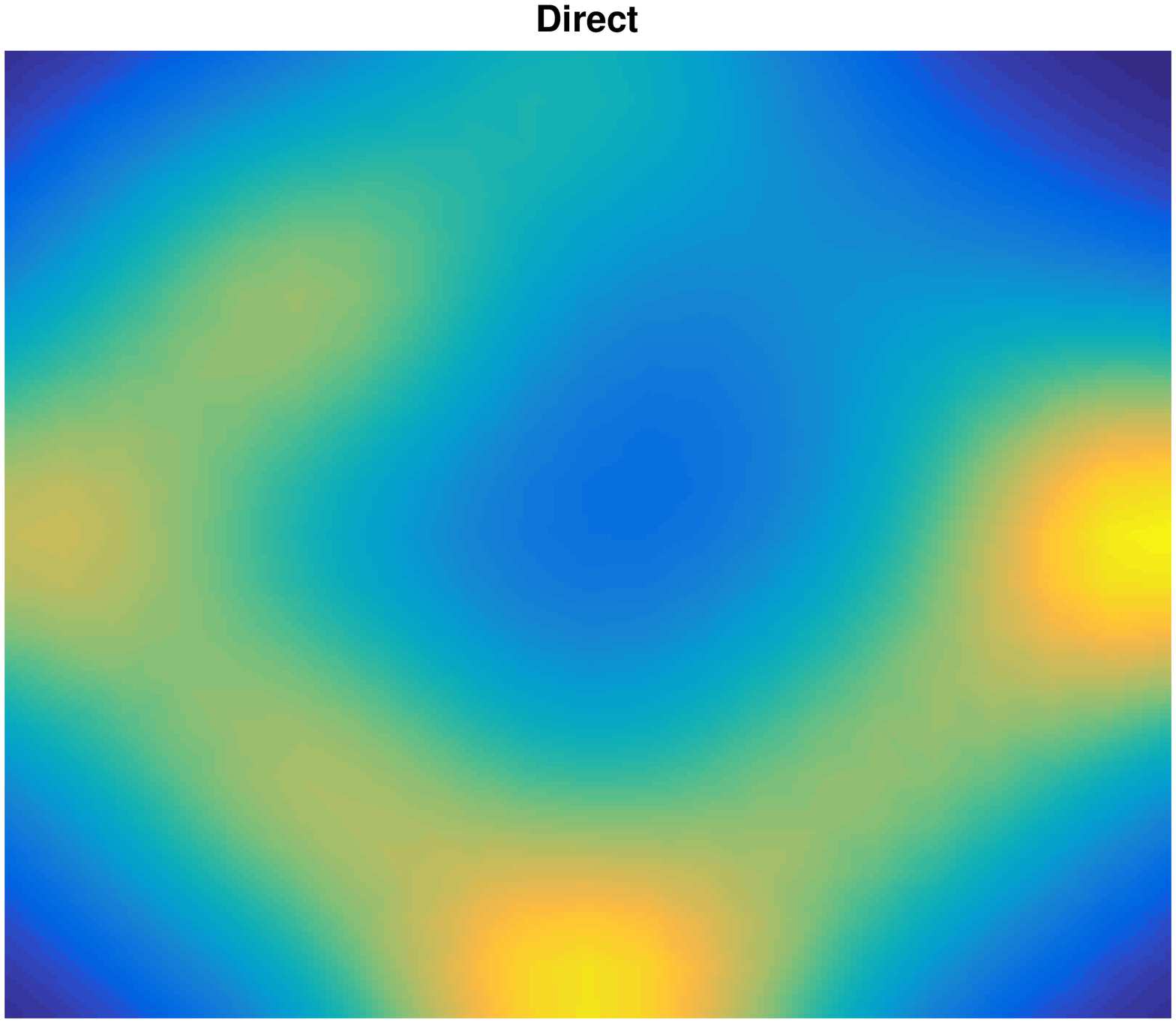}
  \end{overpic}
    \begin{overpic}[width=0.24\textwidth,trim= 20 0 20 15, clip=true,tics=10]{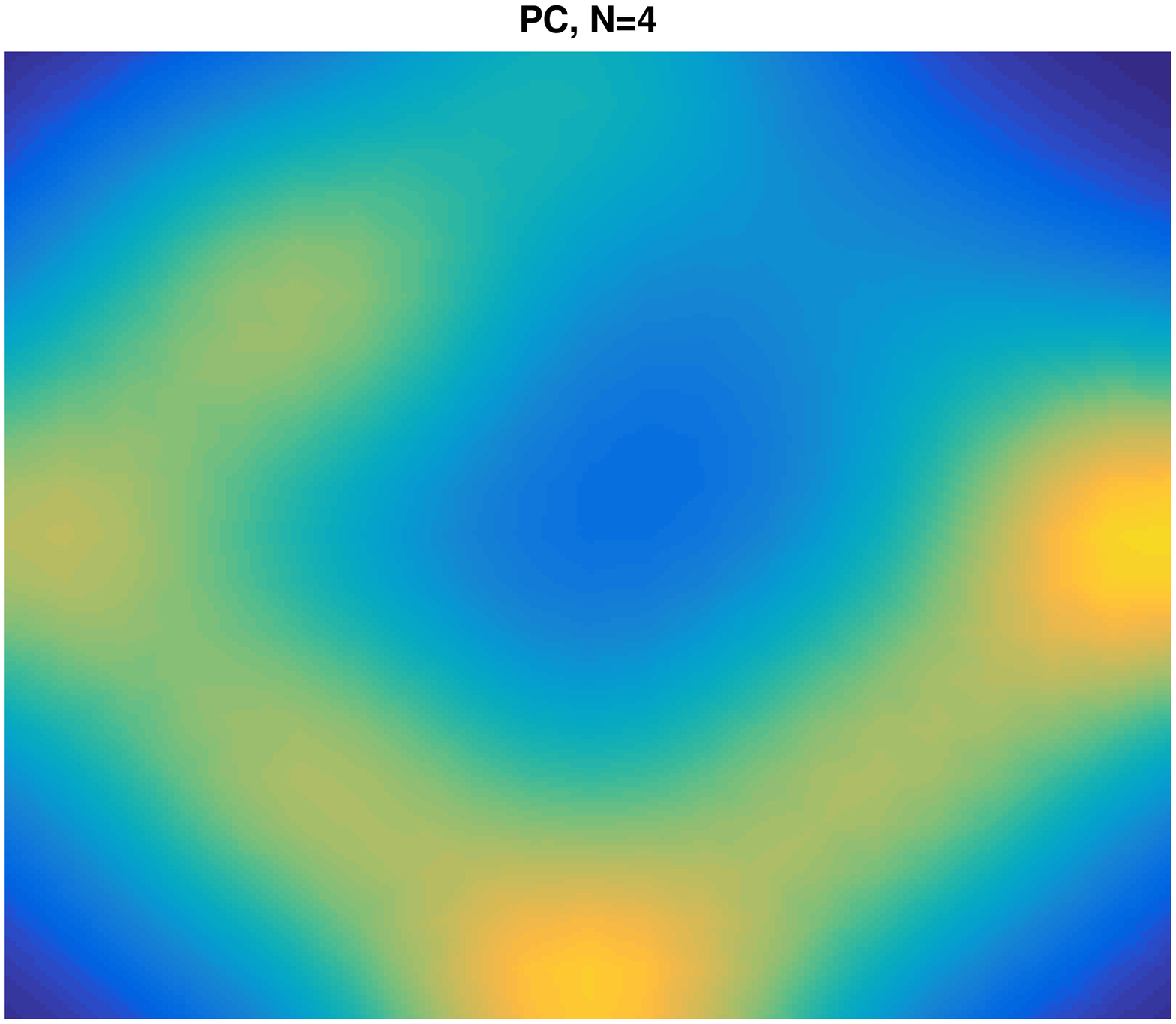}
  \end{overpic}
    \begin{overpic}[width=0.24\textwidth,trim= 20 0 20 15, clip=true,tics=10]{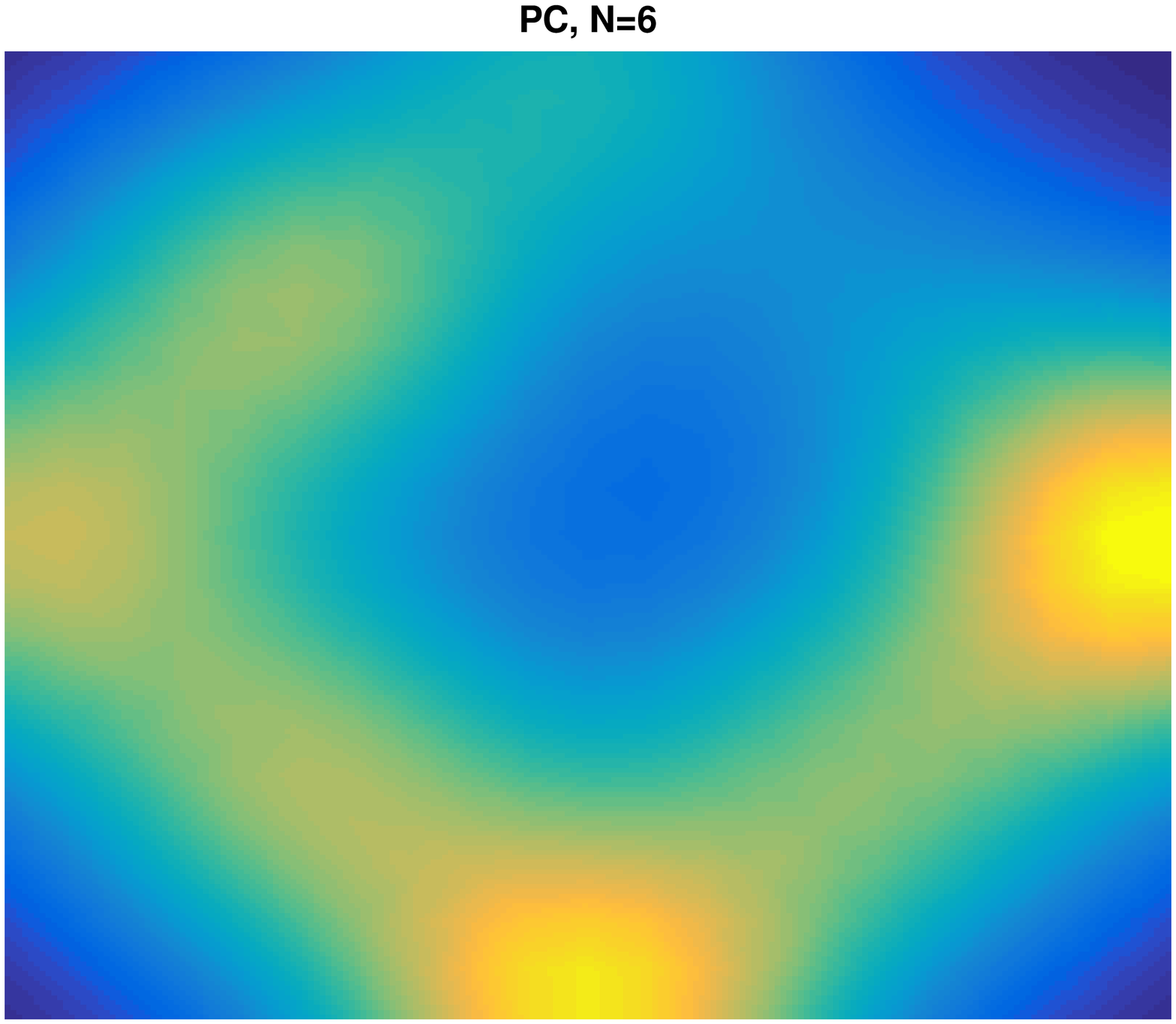}
  \end{overpic}
    \begin{overpic}[width=0.24\textwidth,trim= 20 0 20 15, clip=true,tics=10]{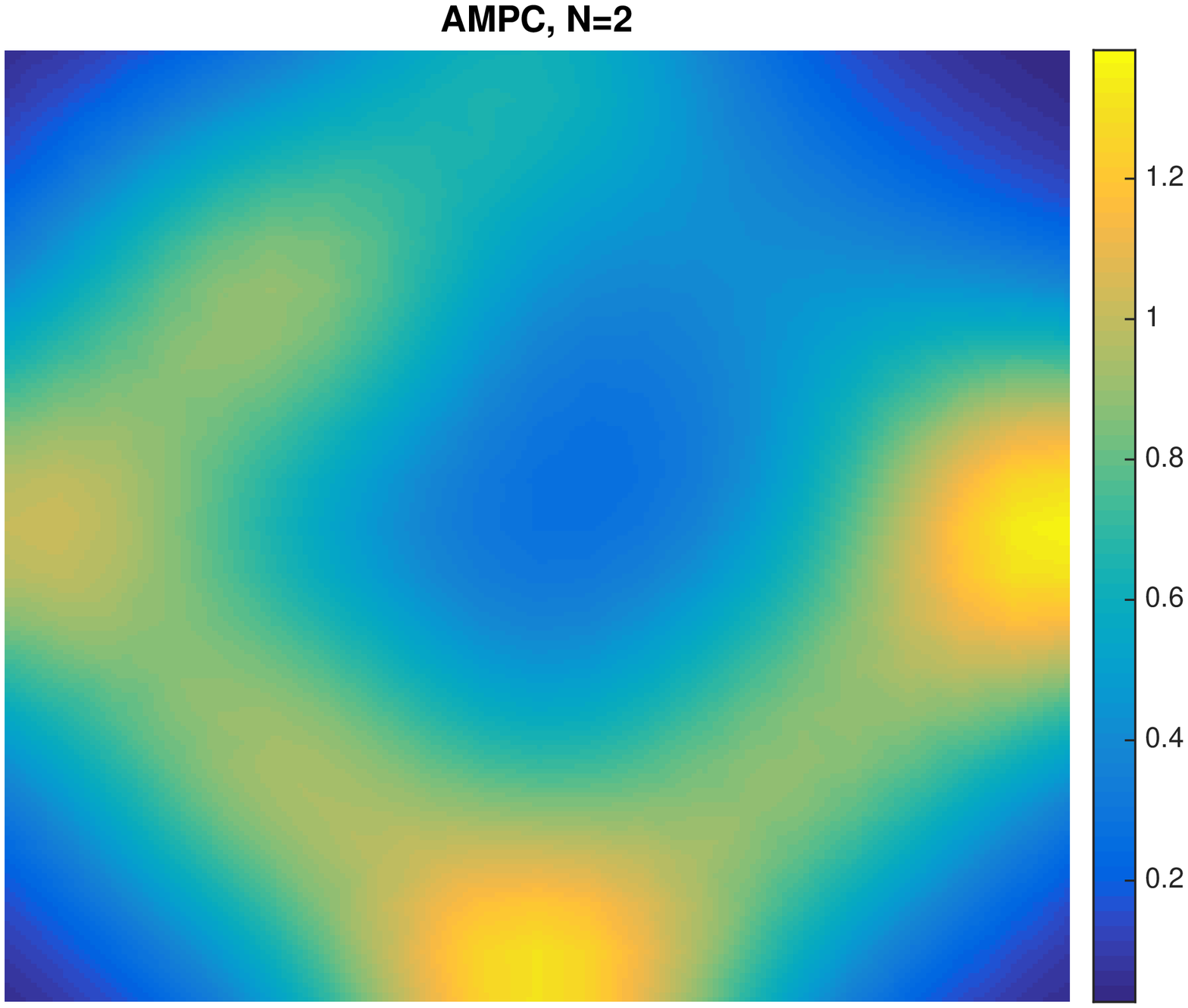}
  \end{overpic}
    \begin{overpic}[width=0.24\textwidth,trim=20 0 20 15, clip=true,tics=10]{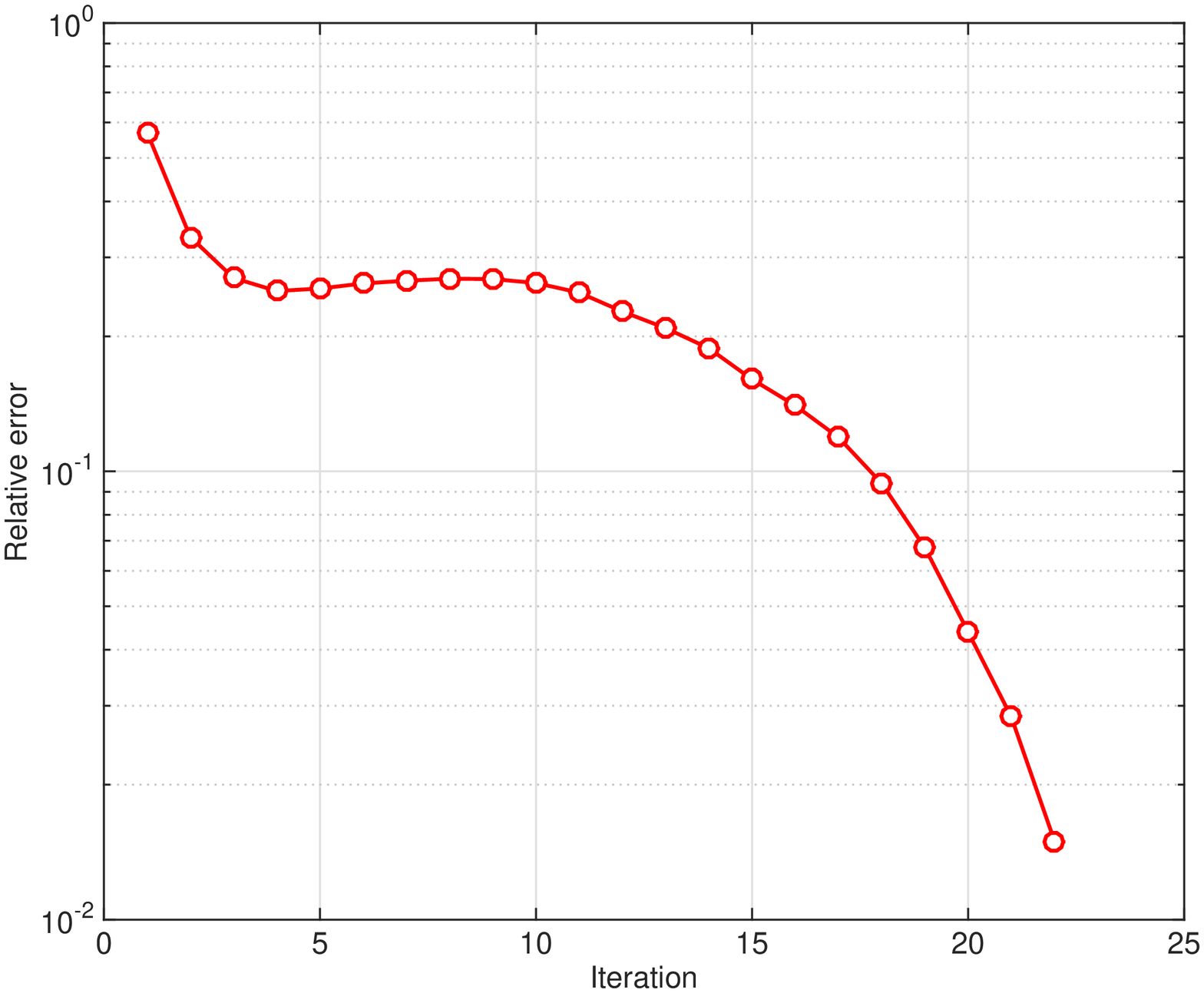}
  \end{overpic}
    \begin{overpic}[width=0.24\textwidth,trim= 20 0 20 15, clip=true,tics=10]{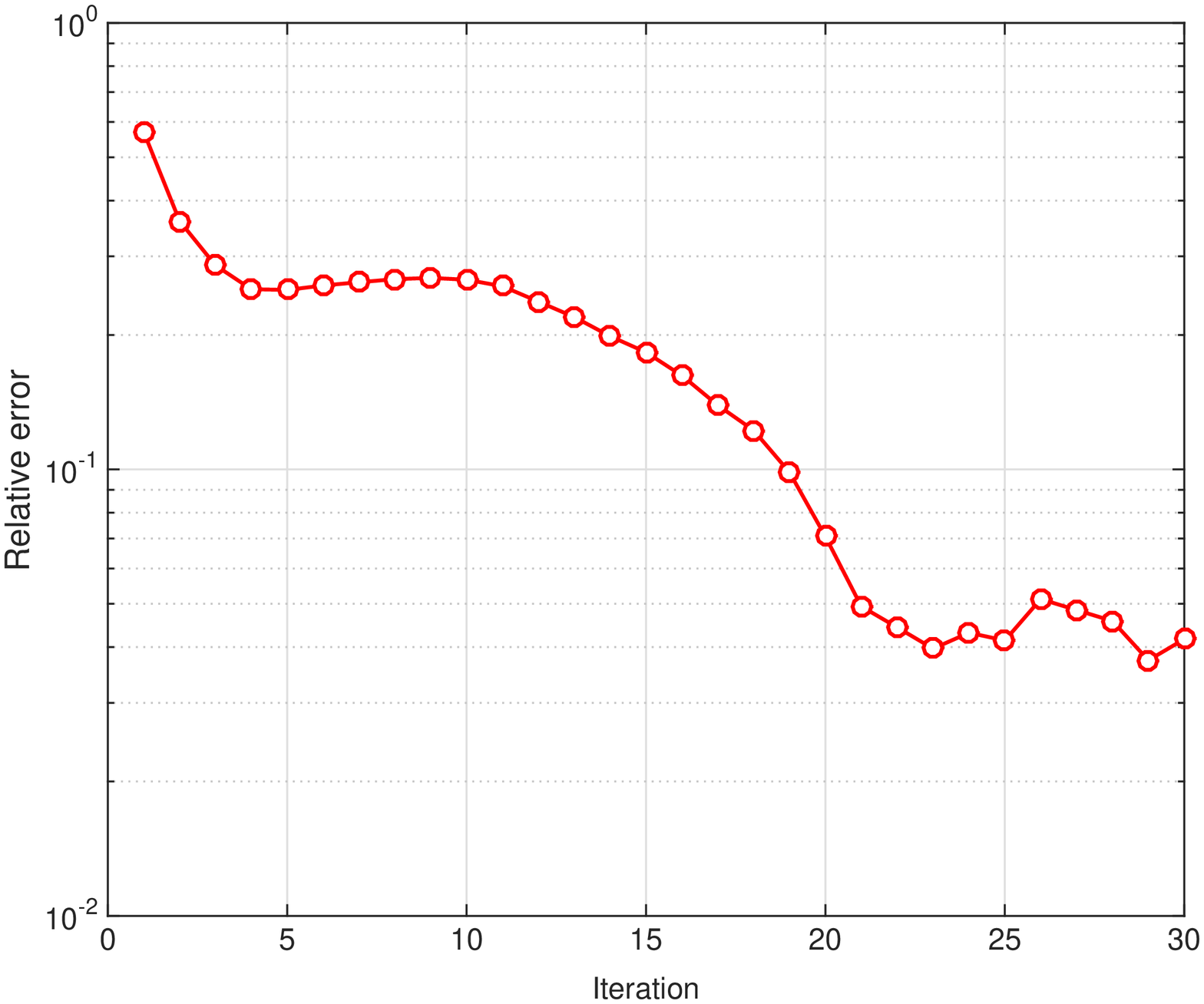}
  \end{overpic}
    \begin{overpic}[width=0.24\textwidth,trim= 20 0 20 15, clip=true,tics=10]{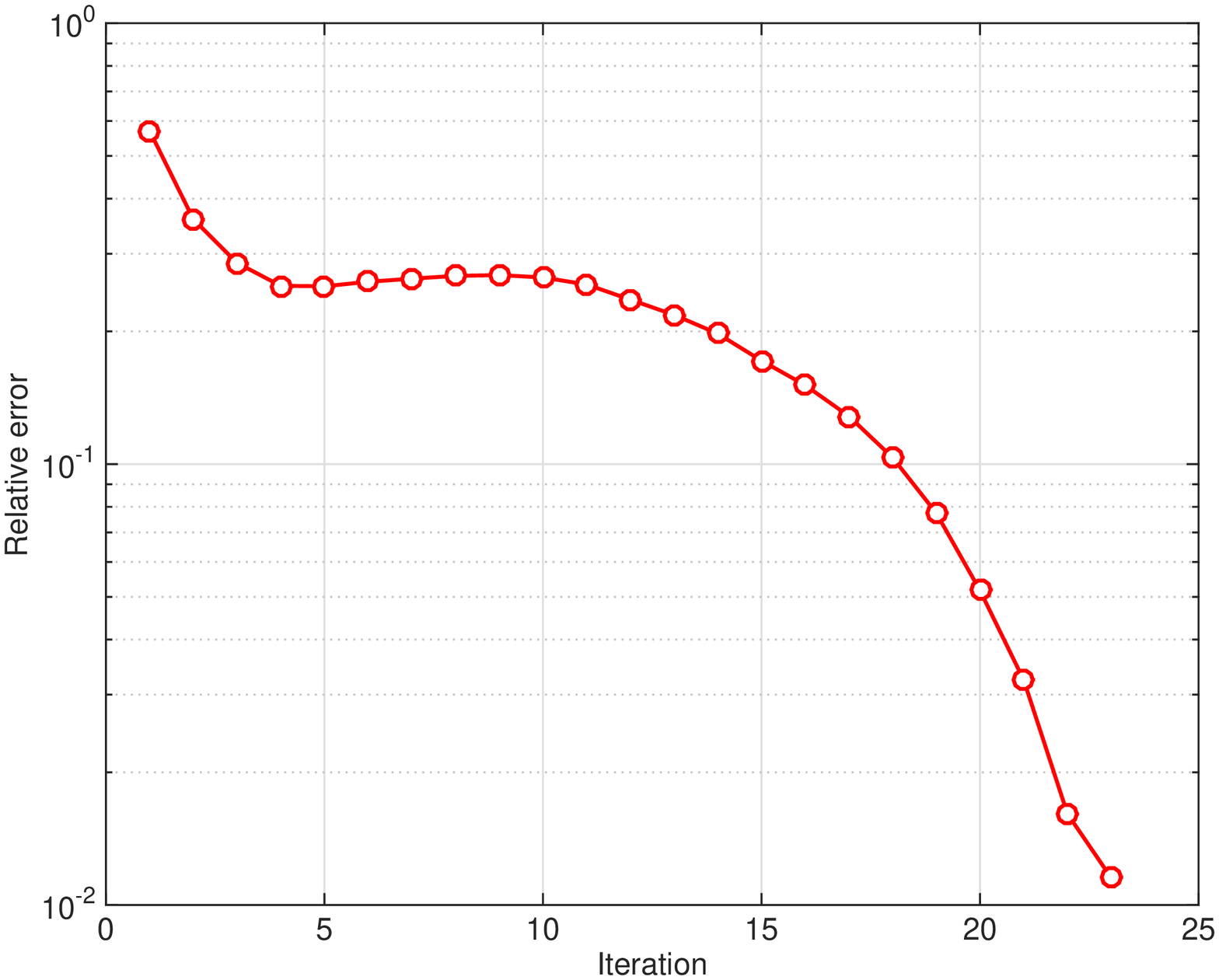}
  \end{overpic}
    \begin{overpic}[width=0.24\textwidth,trim= 20 0 20 15, clip=true,tics=10]{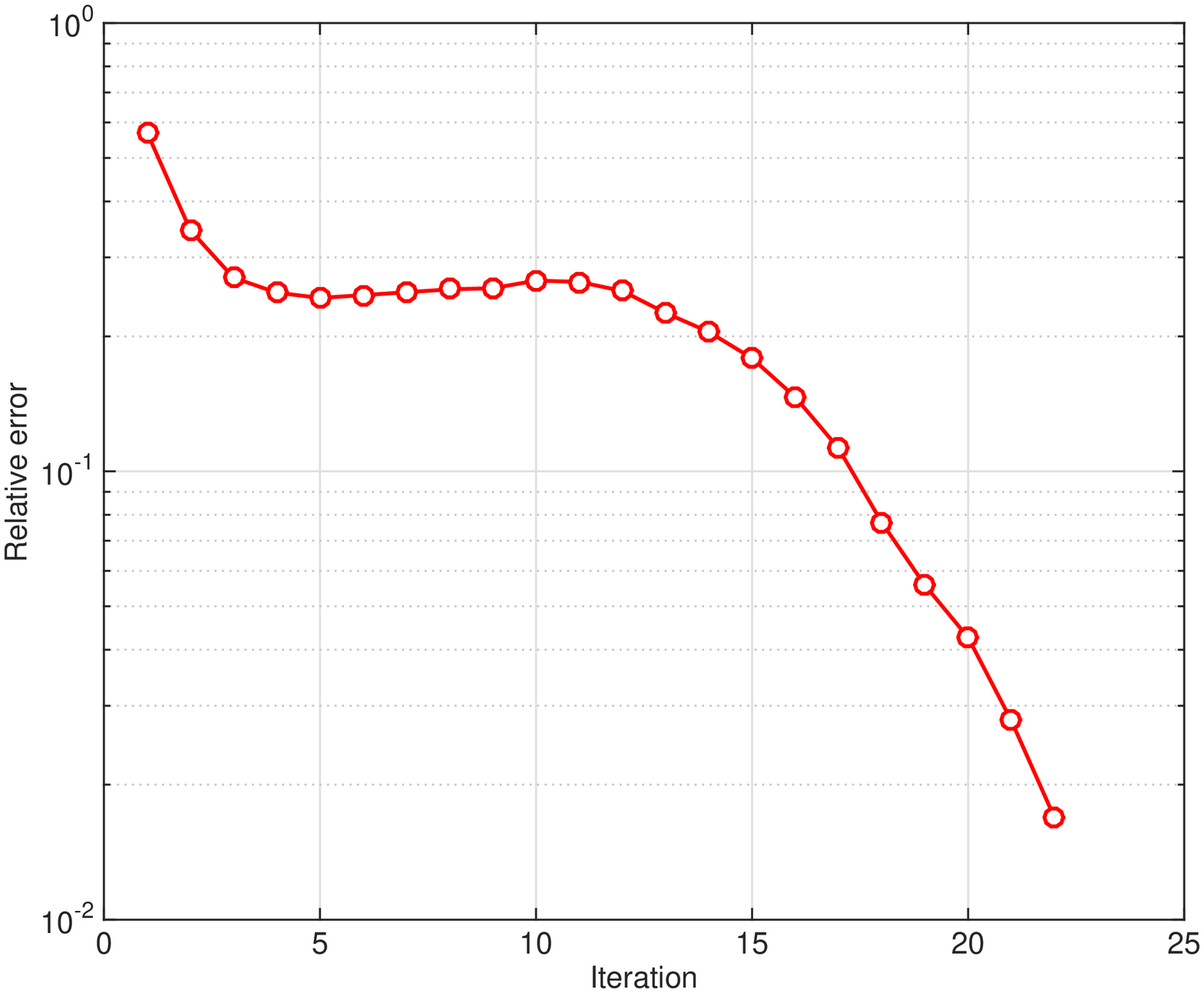}
  \end{overpic}
\end{center}
\caption{Example 2.   Numerical results for the final iteration using $N_e =100$ and different methods (from left to right): Direct; PC(N=4); PC (N=6); AMPC (N=2, $tol=1\times 10^{-3}$).}\label{mean-eg2}
\end{figure}

\begin{figure}
\begin{center}
  \begin{overpic}[width=.32\textwidth,trim=20 0 20 15, clip=true,tics=10]{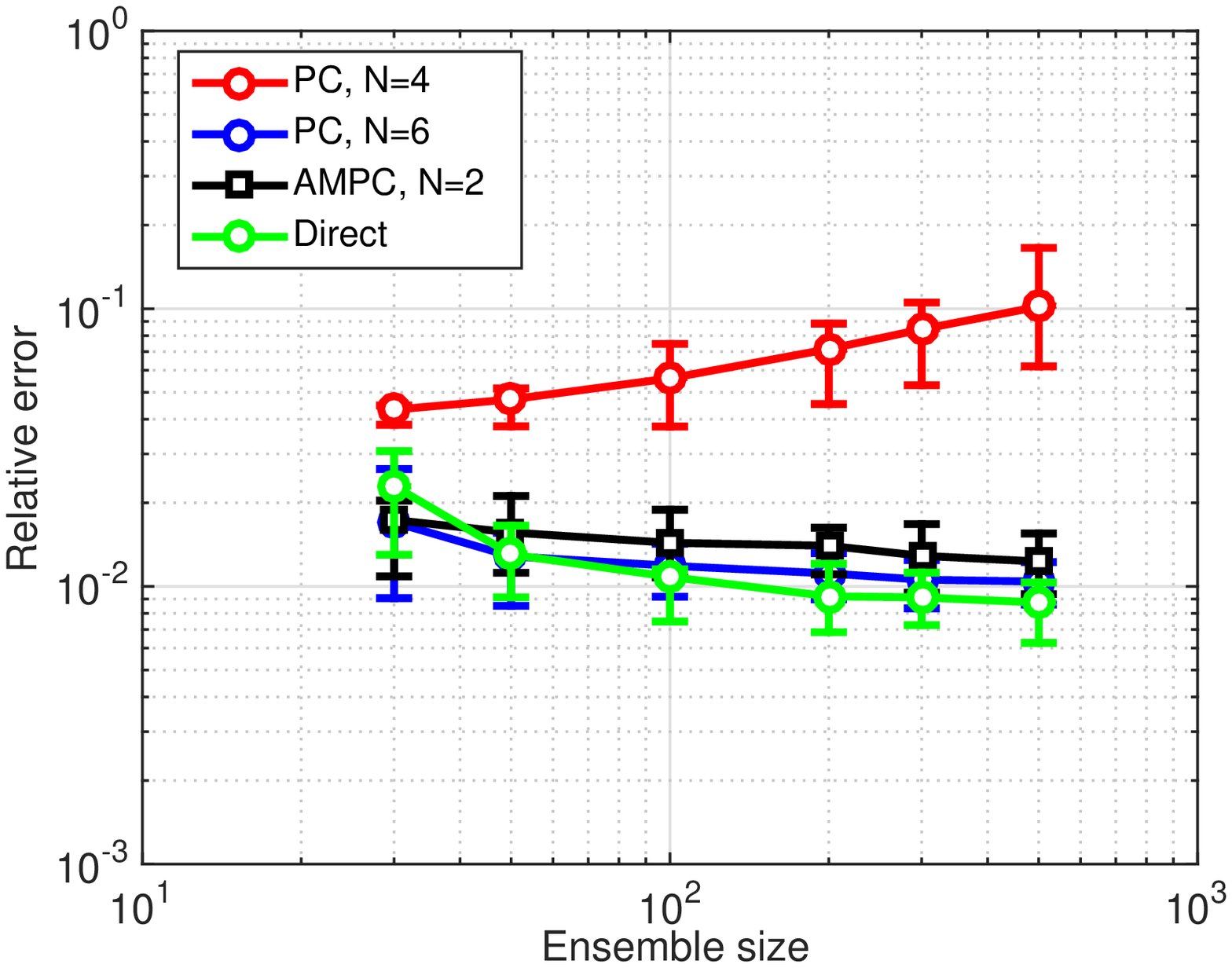}
  \end{overpic}
     \begin{overpic}[width=.32\textwidth,trim=20 0 20 15, clip=true,tics=10]{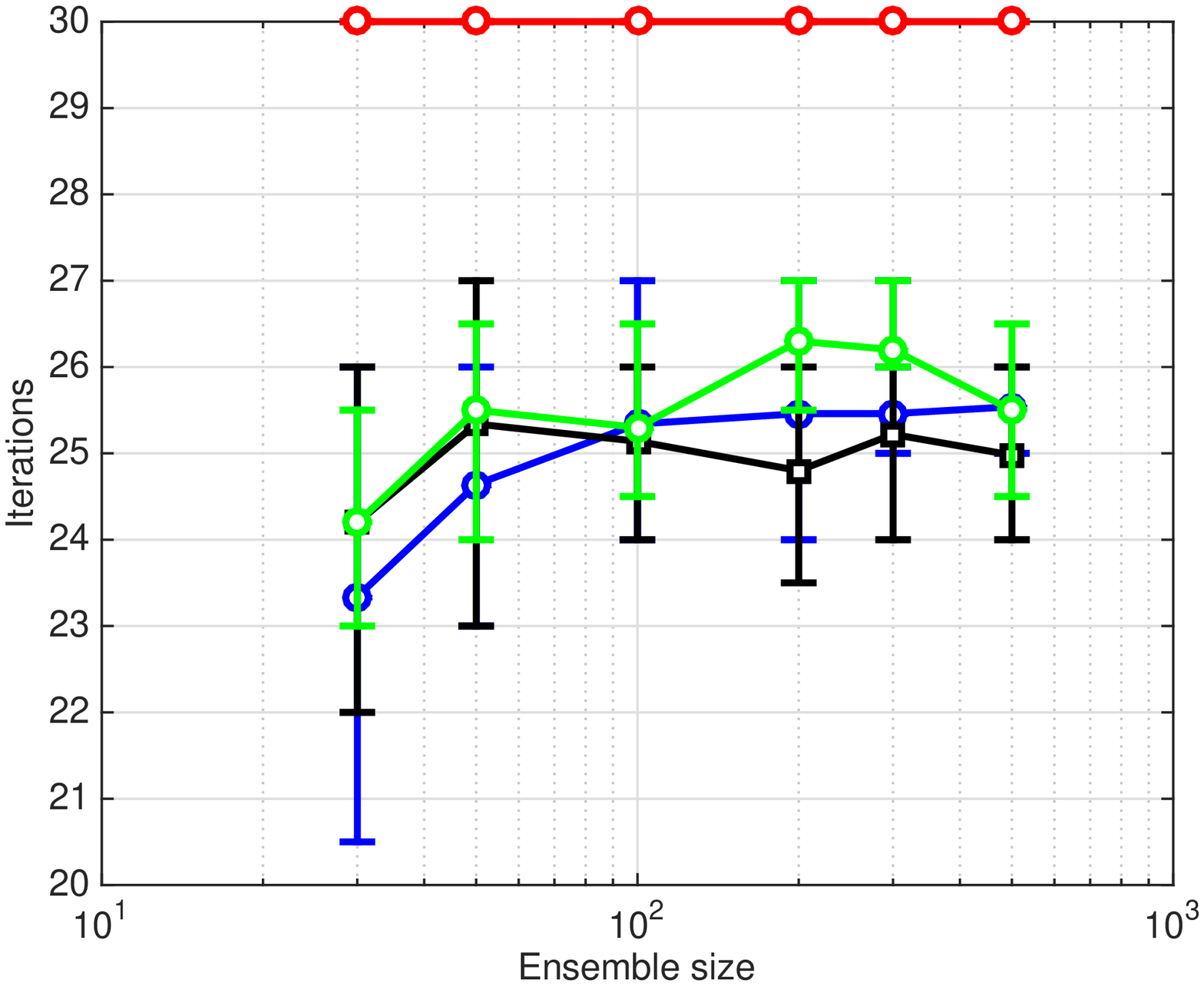}
  \end{overpic}
   \begin{overpic}[width=.32\textwidth,trim=20 0 20 15, clip=true,tics=10]{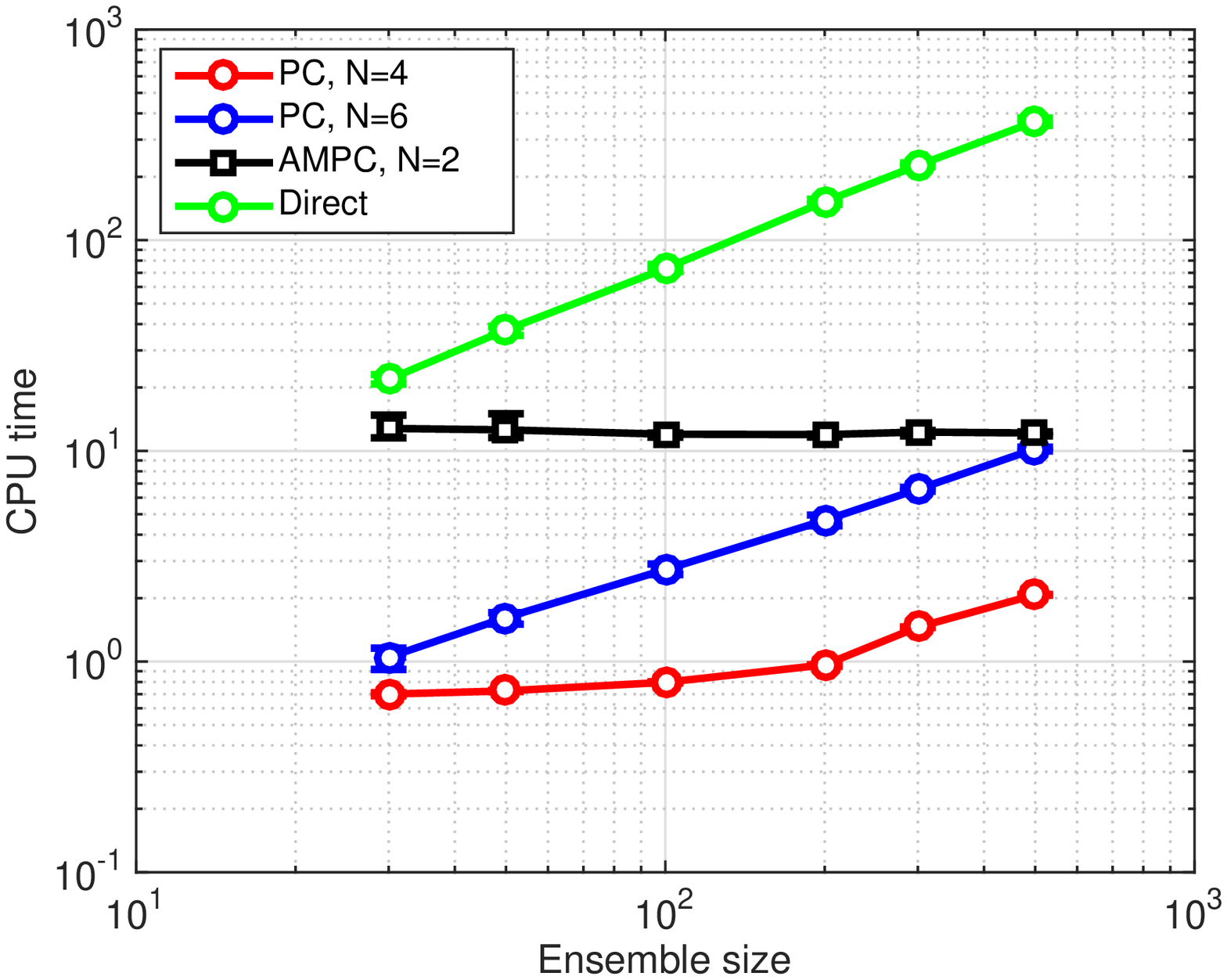}
  \end{overpic}
\end{center}
\caption{Example 2. Numerical results with different ensemble size. Left: values of rel for the estimated results. Middle: number of required iterations. Right: online CPU times.}\label{cpu_eg2}
\end{figure}

Similar to the first example, we numerically investigate the efficiency of the AMPC approach.  Using the same setting as Example 1, we plot the final iteration reconstruction with $N_e =100$ and use the initial mean ensemble given in Figure \ref{exact_eg1}.  The corresponding results are shown in Figures \ref{mean-eg2} and \ref{cpu_eg2}.  Compare with Figure \ref{mean-eg1-PC}, it can be seen that the numerical results obtained by the three approaches are practically identical in this test case,  but the online computing time required by AMPC and PC is only a small fraction of that by the conventional EKI, see the right of Figure \ref{cpu_eg2}.  However, consider the computational cost of the building the PC model,  the total CPU times of AMPC is much smaller than PC-based EKI. This also confirms the efficiency of the AMPC algorithm for this best-case-scenario.

\subsection{Example 3: a high dimensional inverse problem}
\begin{figure}
\begin{center}
  \begin{overpic}[width=.45\textwidth,trim=20 0 20 15, clip=true,tics=10]{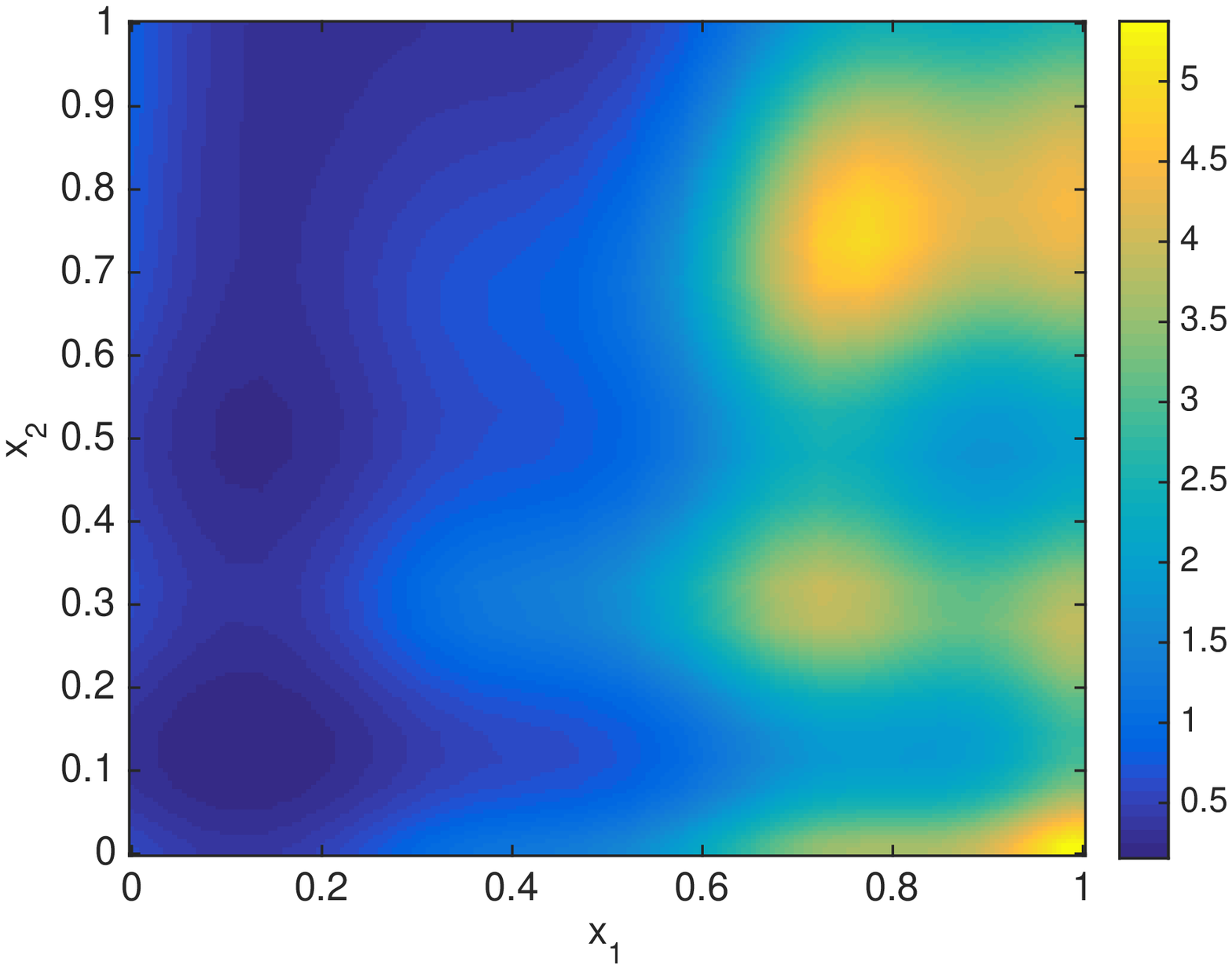}
  \end{overpic}
   \begin{overpic}[width=.45\textwidth,trim=20 0 20 15, clip=true,tics=10]{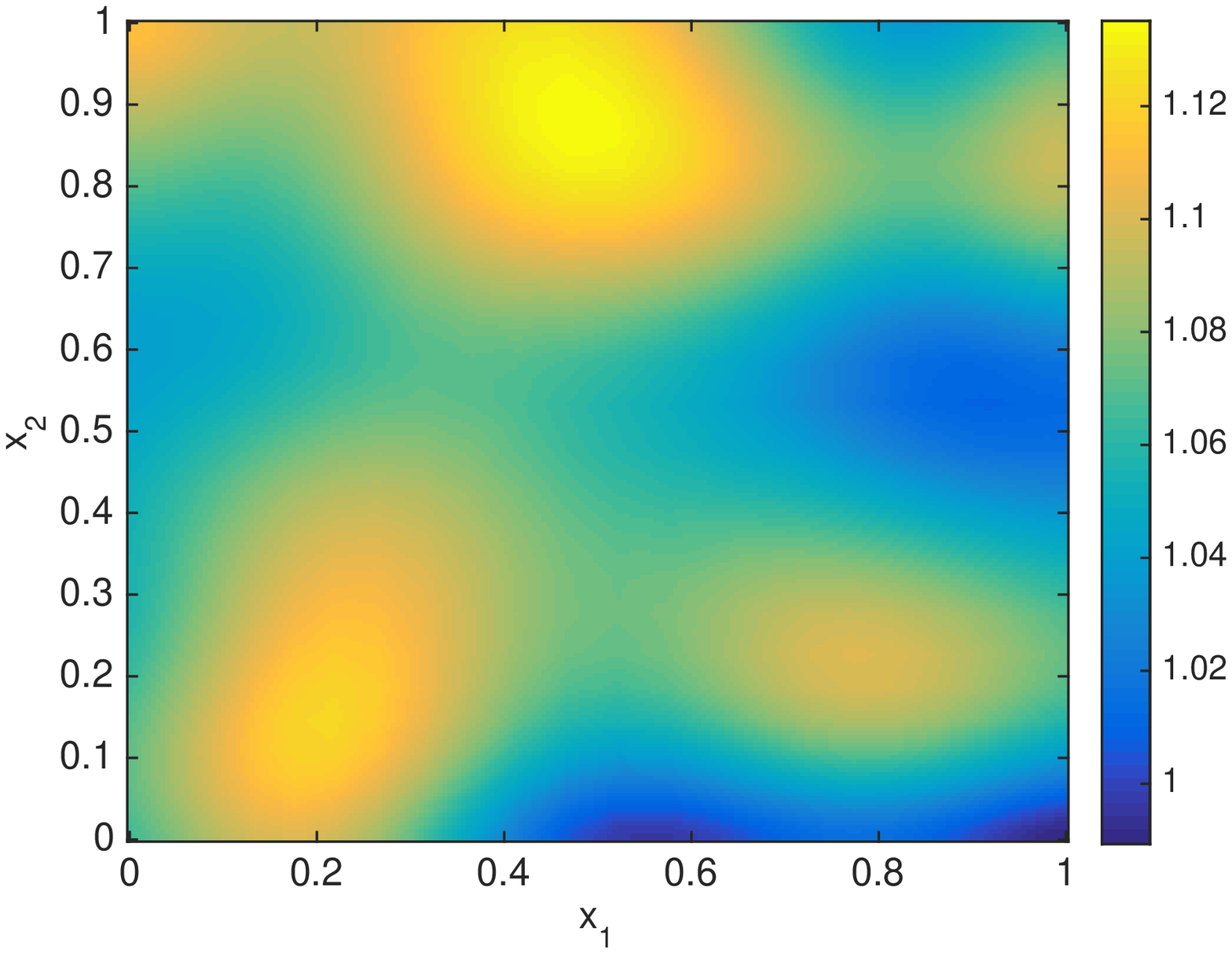}
  \end{overpic}
\end{center}
\caption{Example 3. Left: the true permeability used for generating the synthetic data sets. Right: the initial ensemble mean.}\label{exact_eg3}
\end{figure}

In the last example, we consider the permeabilities as a random field.  Especially, the log-diffusivity field $\log\kappa(x)$ is endowed with a Gaussian process prior, with mean zero and an isotropic exponential covariance kernel:
\begin{equation*}
C(x_1,x_2)=\sigma^2 \exp\Big(-\frac{\|x_1-x_2\|}{2l^2}\Big),
\end{equation*}
for which we choose variance $\sigma^2=1$ and a length scale $l^2 = 0.25$. This prior allows the field to be easily parameterized with a Karhunen-Loeve expansion:
\begin{equation}
\log\kappa(x;\theta) \approx \sum^{d}_{i=1} \theta^i \sqrt{\lambda_i} \phi_i(x),
\end{equation}
where $\lambda_i$ and $\phi_i(x)$ are  the eigenvalues and eigenfunctions, respectively, of the integral operator on $[0,1]^2$ defined by the kernel $C$, and the parameter $\theta^i$ are endowed with independent standard normal priors, $\theta^i \sim N(0,1)$. These parameters then become the targets of inference. In particular, we truncate the Karhunen-Loeve expansion at $d=22$ modes that preserve $95\%$ energy of the prior distribution.  In the numerical simulation, we use the true permeability field that is directly drawn from the prior distribution.   The true permeability field  used to generate the test data, and the initial ensemble mean of the EKI are shown in Figure \ref{exact_eg3}.  The measurement sensors of $u$ are evenly distributed over $\Omega$ with grid spacing 0.1.  Similar to example 1, at each sensor location, three measurements are taken at  time $t=\{0.25,0.75,1\}$. The observational errors are taken to be additive and Gaussian:
\begin{equation*}
y_j = u(x_j,t_j;\theta) +\xi_j,
\end{equation*}
with $\xi_j \sim N(0,0.01^2)$.

 \begin{figure}
\begin{center}
  \begin{overpic}[width=0.24\textwidth,trim=20 0 20 15, clip=true,tics=10]{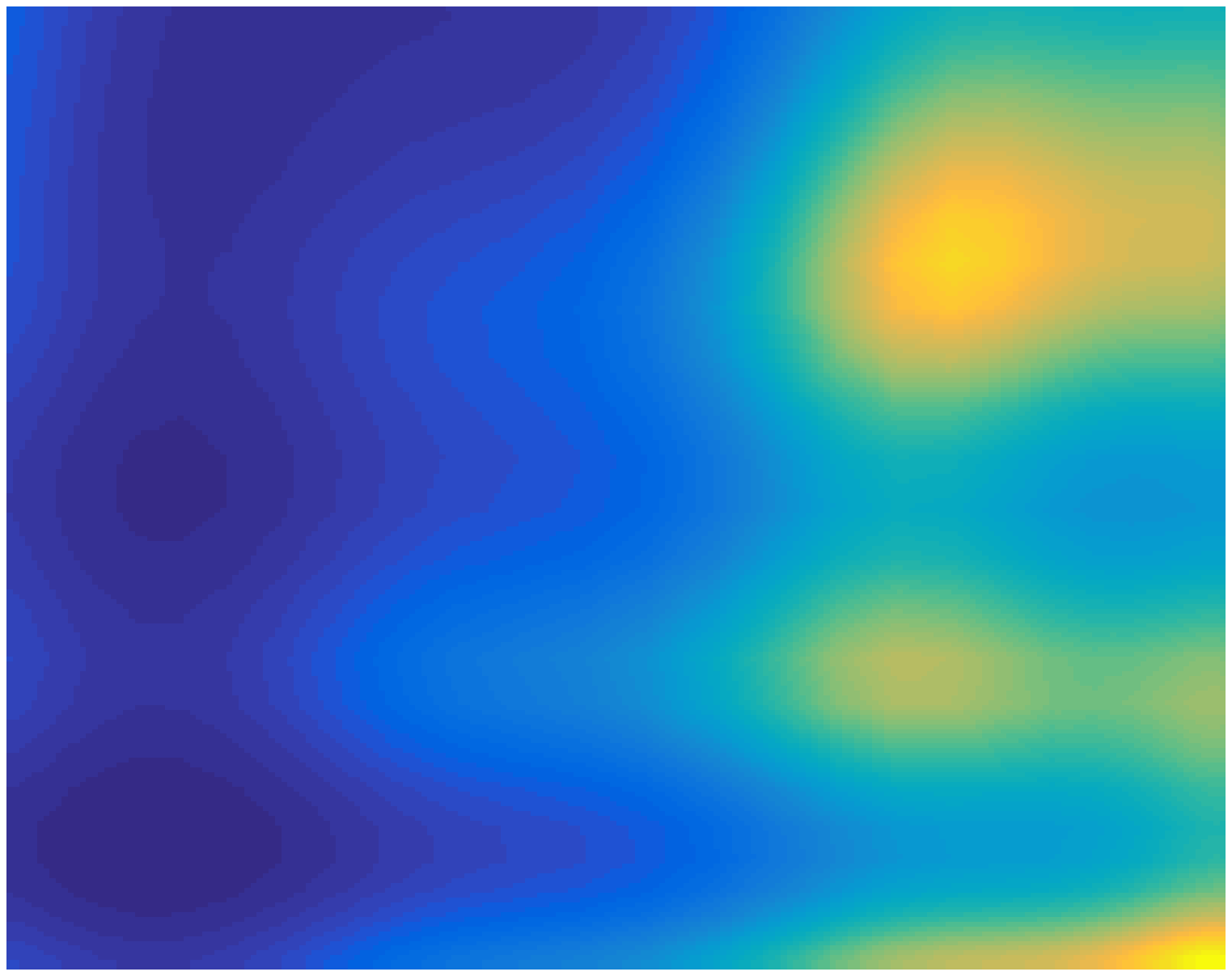}
  \end{overpic}
    \begin{overpic}[width=0.24\textwidth,trim= 20 0 20 15, clip=true,tics=10]{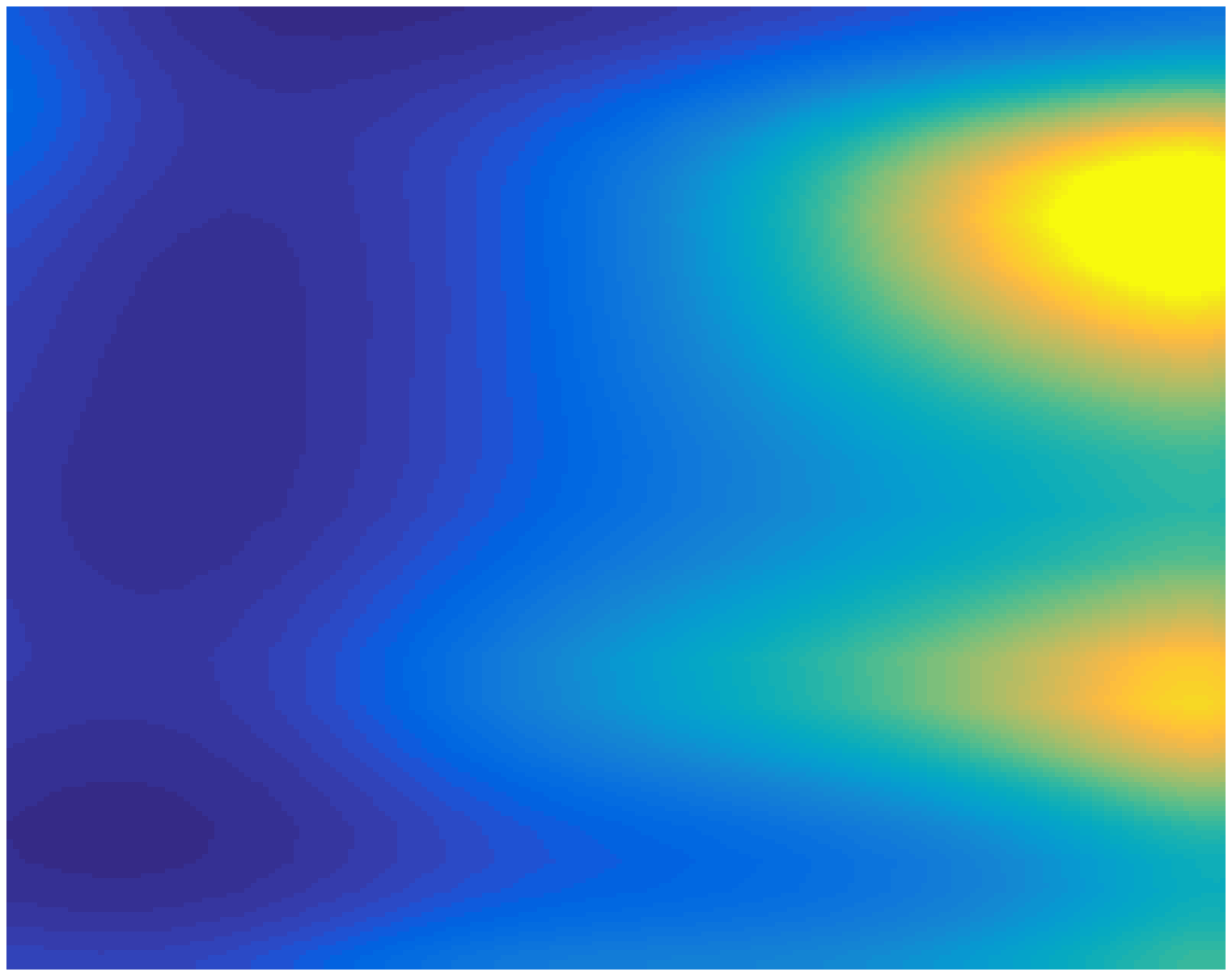}
  \end{overpic}
    \begin{overpic}[width=0.24\textwidth,trim= 20 0 20 15, clip=true,tics=10]{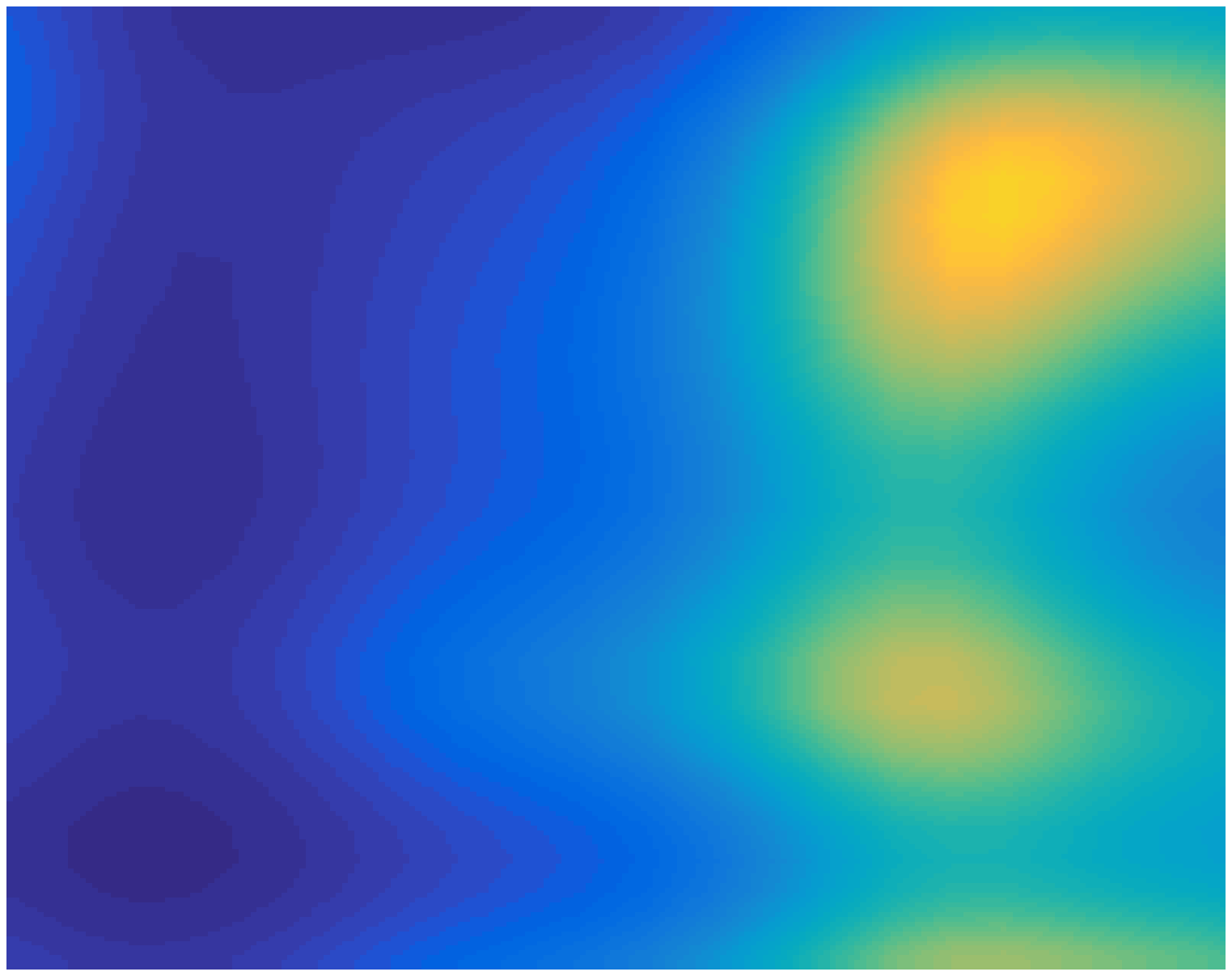}
  \end{overpic}
    \begin{overpic}[width=0.24\textwidth,trim= 20 0 20 15, clip=true,tics=10]{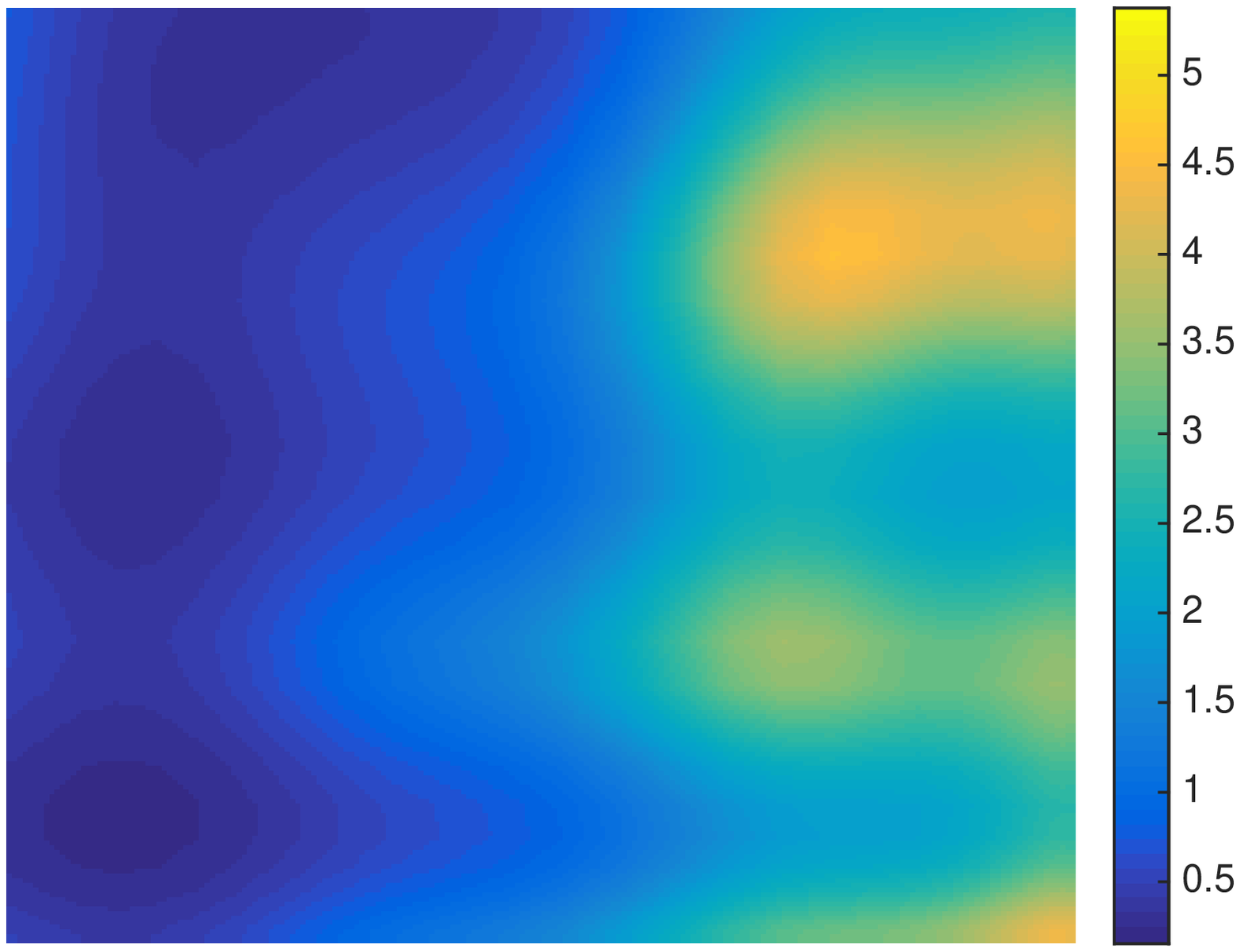}
  \end{overpic}
    \begin{overpic}[width=0.24\textwidth,trim=20 0 20 15, clip=true,tics=10]{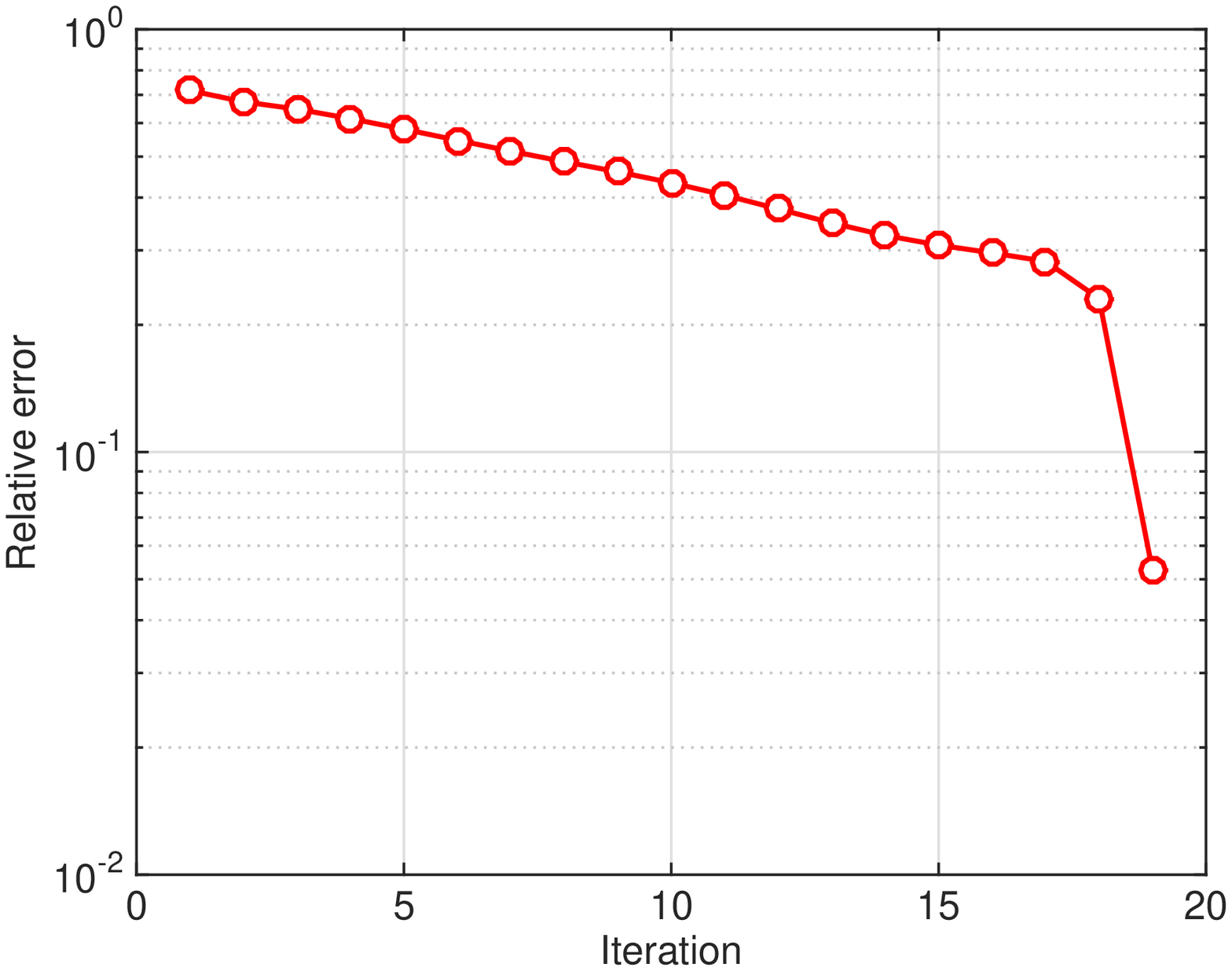}
  \end{overpic}
    \begin{overpic}[width=0.24\textwidth,trim= 20 0 20 15, clip=true,tics=10]{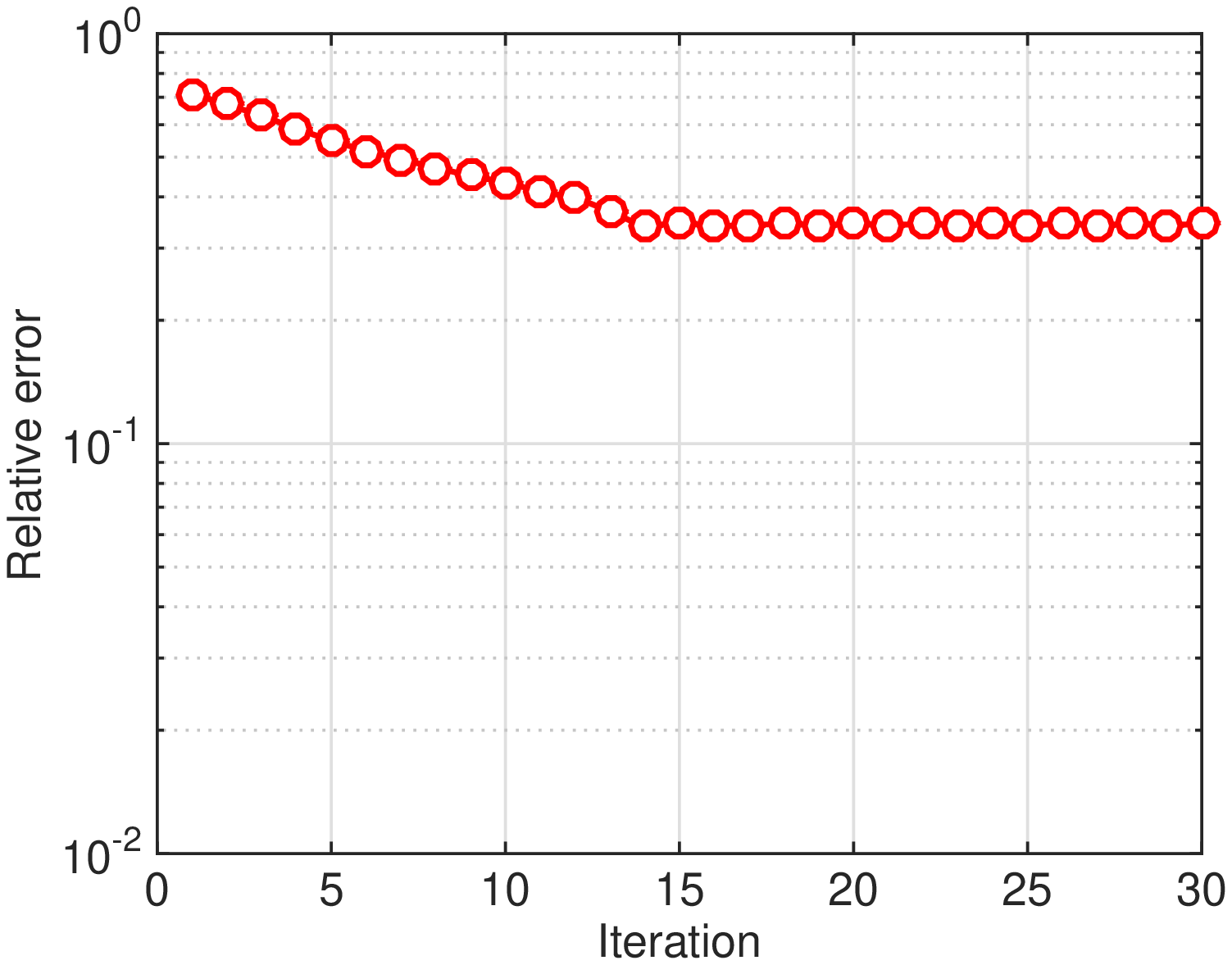}
  \end{overpic}
    \begin{overpic}[width=0.24\textwidth,trim= 20 0 20 15, clip=true,tics=10]{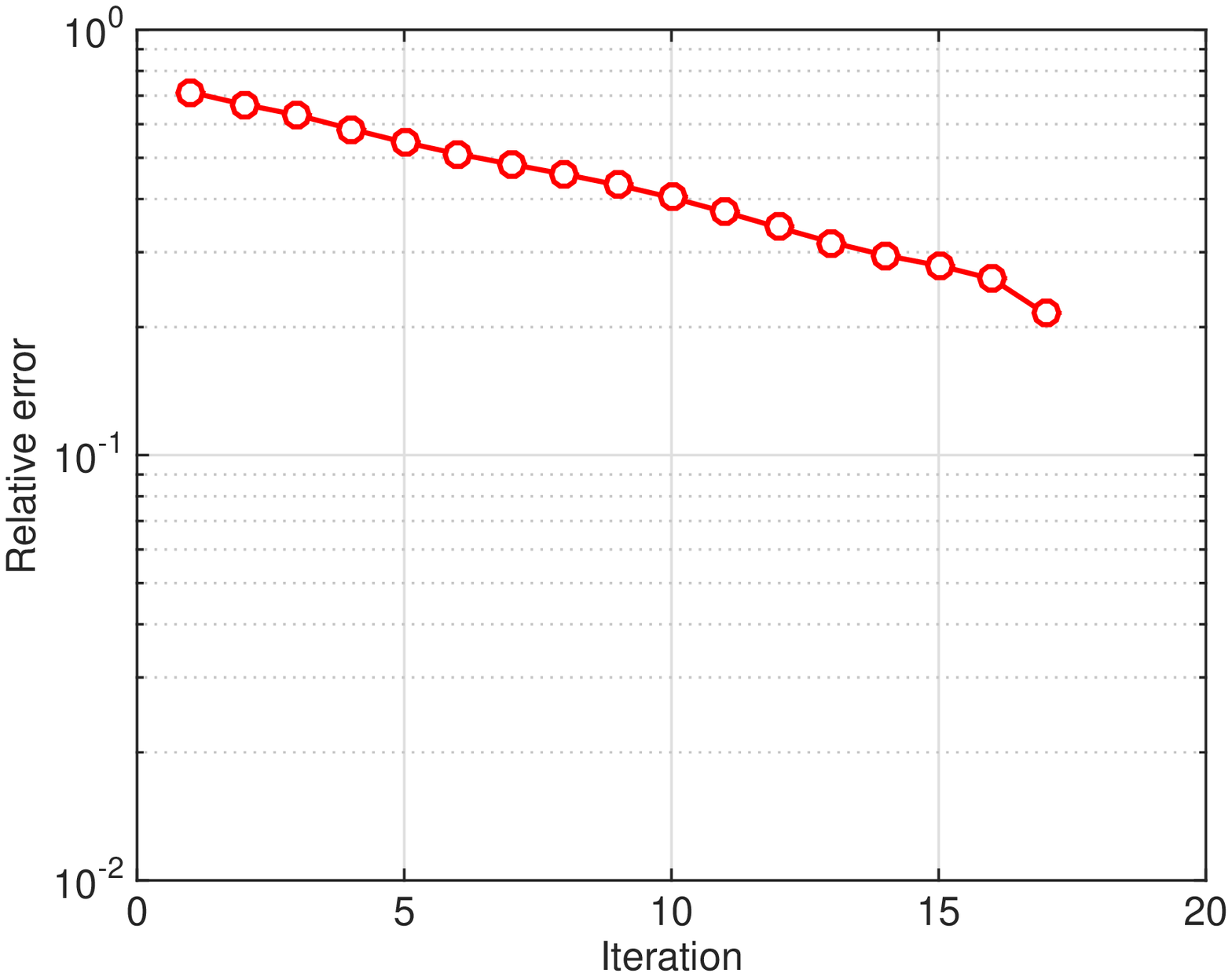}
  \end{overpic}
    \begin{overpic}[width=0.24\textwidth,trim= 20 0 20 15, clip=true,tics=10]{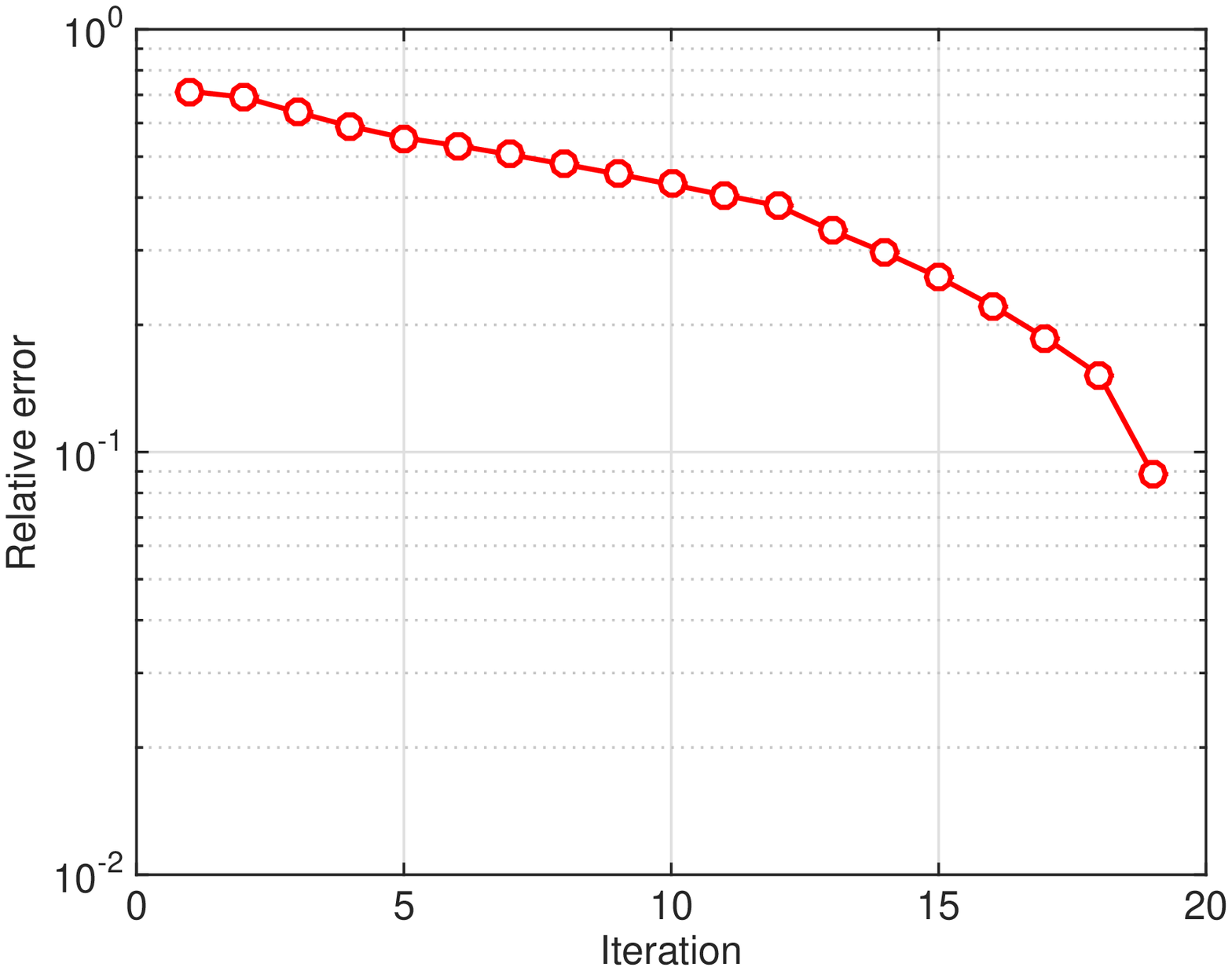}
  \end{overpic}
\end{center}

\caption{Example 3.   Numerical results for the final iteration using $N_e =300$ and different methods(from left to right): Direct; PC (N=2); PC (N=3); AMPC (N=1, $tol=1\times 10^{-2}$).}\label{mean-eg3}
\end{figure}

Figure \ref{mean-eg3} plots the final iteration reconstruction with $N_e =300$ using three different approaches.  As expected, a poor estimate is obtained by the PC-based EKI with a lower order $N=2$ ($rel=0.3430$). The results are improved with $N=3$ ($rel =0.2146$).  In contrast, even with a lower order $N=1$, the numerical  results obtained by AMPC agree well with the exact solution $(rel =0.0889)$.

The computational costs and the relative errors $rel$ of the final iteration  for the different algorithms are shown in Table \ref{eg3_time}.  Building a PC surrogate of order $N=2$ (resp. N=3) requires an offline CPU time of 29.73 (resp. 256.54), whereas its online evaluation requires 2.99s (resp. 6.66s). This fact has a major drawback for  PC-based EKI to solve high-dimensional problems: the total CPU times increase fast with respect to the polynomial order $N$.    On the other hand, for the AMPC algorithm with the PC order $N=1$ and $tol=1\times 10^{-2}$, the  offline and online CPU times are 2.69s and 10.35s, respectively. And the  relative error $rel$ of AMPC is about $0.0889$, which is more efficient than the PC-based EKI with order $N=3$.  This demonstrated that the AMPC can provide with much more accurate results, yet with less computational time.

 \begin{table}[tp]
      \caption{Example 3. Computational times, in seconds, given by three different methods. $tol=1\times10^{-2}, N_e=300$. }\label{eg3_time}
  \centering
  \fontsize{6}{12}\selectfont
  \begin{threeparttable}
    \begin{tabular}{ c cccccc}
  \toprule
 & \multicolumn{2}{c}{Offline}&\multicolumn{2}{c}{Online}\cr
\cmidrule(lr){2-3} \cmidrule(lr){4-5}

  \multirow{1}{*}{Method}  &$\text{$\#$ of model eval.}$&CPU(s) &$\text{$\#$ of model eval.}$&CPU(s)     &\multirow{1}{*}{Total time(s)}&\multirow{1}{*}{rel}\cr
  \midrule
    Direct                                     & $-$       & $-$         & 5700      &309.19      & 309.19  & 0.0523\cr
   PC, $N=3$                             & 4600  & 256.54     & $-$          &6.66          & 263.20  &0.2146\cr
   PC, $N=2$                             & 552   & 29.73        & $-$          &2.99          & 32.72  & 0.3430\cr
   AMPC, $N= N_C=1$          & 46     & 2.69          & 157          &10.35       & 13.04  &0.0889\cr
      \bottomrule
      \end{tabular}
    \end{threeparttable}

\end{table}

\section{Summary} \label{sec:summary}
In this paper, we developed an adaptive  multi-fidelity PC based EKI algorithm to  solve nonlinear inverse problems. This new strategy combines a large number of low-order PC surrogate model evaluations and a small number of  forward model evaluations, yielding a multi-fidelity approach. The key idea is to construct and refine the multi-fidelity PC surrogate using the updated parameters at each iteration. Then the prediction steps of the EKI are calculated from a large number of realizations generated by the multi-fidelity PC surrogate with virtually no additional  computational cost.  During the entire EKI simulation, the high-fidelity model evaluations are only needed at the refine the multi-fidelity PC model, whose number is much smaller than the total ensemble size of classic EKI. Thus, the computational cost can be significantly reduced.  The performance of the proposed strategy has been illustrated by three numerical examples. Although only the iterative regularizing ensemble Kalman smoother algorithm are considered in this paper, the AMPC scheme  can be conveniently extended to a much wider class of EKIs with simple and minor modifications. The extension of the present algorithm  to highly nonlinear and complex systems is also straightforward.

\bibliographystyle{plain}
\bibliography{mEnPCK}

\end{document}